\documentclass{amsart}
\numberwithin{equation}{section}
\usepackage{psfrag,comment,amssymb,amscd}
\usepackage[all]{xy}
\let\cal\mathcal
\def\Ascr{{\cal A}}
\def\Bscr{{\cal B}}
\def\Cscr{{\cal C}}
\def\Dscr{{\cal D}}

\def\Fscr{{\cal F}}

\def\Lscr{{\cal L}}

\def\Oscr{{\cal O}}

\def\Tscr{{\cal T}}
\def\Uscr{{\cal U}}

\def\id{\text{id}}
\def\Id{\operatorname{id}}

\def\Der{\operatorname{Der}}

\def\ctimes{\mathbin{\hat{\otimes}}}

\def\gr{\operatorname{gr}}
\def\Lie{\mathop{\text{Lie}}}

\def\gr{\operatorname {gr}}
\def\Spec{\operatorname {Spec}}

\def\GL{\operatorname {GL}}

\def\Ext{\operatorname {Ext}}
\def\Hom{\operatorname {Hom}}
\def\uHom{\operatorname {\mathcal{H}\mathit{om}}}

\def\im{\operatorname {im}}

\def\ker{\operatorname {ker}}

\def\id{{\operatorname {id}}}

\def\r{\rightarrow}

\let\invlim\projlim

\newtheorem{lemma}{Lemma}[section]
\newtheorem{proposition}[lemma]{Proposition}
\newtheorem{theorem}[lemma]{Theorem}

\newtheorem{lemmas}{Lemma}[subsection]
\newtheorem{propositions}[lemmas]{Proposition}
\newtheorem{theorems}[lemmas]{Theorem}
\newtheorem{corollarys}[lemmas]{Corollary}

\newtheorem*{sublemma}{Sublemma}

\theoremstyle{definition}

\newtheorem{example}[lemma]{Example}

\newtheorem{examples}[lemmas]{Example}

{

}

\theoremstyle{remark}

\newtheorem{remark}[lemma]{Remark}
\newtheorem{remarks}[lemmas]{Remark}

\newdimen\uboxsep \uboxsep=1ex
\def\uboxn#1{\vtop to 0pt{\hrule height 0pt depth 0pt\vskip\uboxsep
\hbox to 0pt{\hss #1\hss}\vss}}

\def\uboxs#1{\vbox to 0pt{\vss\hbox to 0pt{\hss #1\hss}
\vskip\uboxsep\hrule height 0pt depth 0pt}}

\def\poly{\operatorname{poly}}

\def\coord{\operatorname{coord}}
\def\aff{\operatorname{aff}}
\def\afff{\mathfrak{a}}

\def\poly{\operatorname{poly}}

\def\Aff{\operatorname{\Ascr}}
\def\poly{\operatorname{poly}}
\def\coder{\operatorname{coder}}
\def\cup{\operatorname{cup}}
\def\Harr{\operatorname{Harr}}
\def\Diff{\operatorname{Diff}}
\def\Poiss{\operatorname{Poiss}}
\def\Lie{\operatorname{Lie}}
\def\HKR{\operatorname{HKR}}

\keywords{Deformation quantization}
\author{Damien Calaque}
\thanks{At the final stage of this work, the first author has been supported by the European Union thanks to 
a Marie Curie Intra-European fellowship (contract number MEIF-CT-2007-042212). }
\address{Department of Mathematics \\ ETH Zurich \\ 8092 Zurich \\ Switzerland}
\address{Universit\'e de Lyon \\ Universit\'e Lyon 1 \\ Institut Camille Jordan CNRS UMR 5208 \\
43 boulevard du 11 novembre 1918 \\ F-69622 Villeurbanne Cedex}
\email{calaque@math.univ-lyon1.fr, damien.calaque@math.ethz.ch}
\author{Michel Van den Bergh}
\address{Departement WNI\\Universiteit Hasselt\\ Universitaire Campus\\ Building D\\ 3590 
Diepenbeek\\ Belgium} 
\thanks{The second author is a director of research at the FWO} 
\thanks{The results of this paper were partially obtained while the
second author was visiting the Universit\'e Claude Bernard at Lyon. He hereby
thanks the latter for its kind hospitality.}
\email{michel.vandenbergh@uhasselt.be} 

\dedicatory{Dedicated to Pierre Deligne's 65-th birthday}
\title{Global formality at the $G_\infty$-level}

\begin{document}

\begin{abstract}
In this paper we prove that the sheaf of $\Lscr$-poly-differential operators
for a locally free Lie algebroid $\Lscr$ is formal when viewed
as a sheaf of $G_\infty$-algebras via Tamarkin's morphism of DG-operads
$G_\infty\rightarrow B_\infty$.

In an appendix  we prove a strengthening of Halbout's globalization
result for Tamarkin's local quasi-isomorphism. 
\end{abstract}

\maketitle

\setcounter{tocdepth}{1}

\tableofcontents

\section{Introduction}

Throughout $k$ is a field of characteristic zero. In this paper we
extend the global $L_\infty$-formality result \cite[Thm
7.4.1]{vdbcalaque} (see also \cite{ye3}) for Lie algebroids to the
$G_\infty$-context.

A similar global $G_\infty$-formality result was obtained in
\cite{dtt} when the Lie algebroid is the tangent bundle of a smooth
space (in a suitable context). The methods in loc.\ cit.\ are however
quite different. In the special case of a $C^\infty$-manifold, 
 the result has also been obtained in \cite{GH}.

Throughout $\Cscr$ is a site equipped with a sheaf of commutative,
associative $k$-algebras $\Oscr$.  Furthermore $\Lscr$ is a locally
free Lie algebroid over $(\Cscr,\Oscr)$ of constant rank~$d$.

Let $T_{\poly}^\Lscr(\Oscr)$, $D_{\poly}^\Lscr(\Oscr)$ be respectively
the sheaves of $\Lscr$-poly-vector fields and
$\Lscr$-poly-differential operators on $\Oscr$ (see \cite{cal} or
\S\ref{ref-8-30} below). Then $T_{\poly}^\Lscr(\Oscr)$ is canonically a sheaf of
Gerstenhaber algebras and $D_{\poly}^\Lscr(\Oscr)$ is canonically a
$B_\infty$-algebra (see \cite{VG} or \S\ref{ref-3.1-2} below).  

Tamarkin constructs  in \cite{Tamarkin} (see \S\ref{ref-5-9})
a remarkable morphism of operads $G_\infty\rightarrow B_\infty$ which is
intimately connected to the celebrated (now multiply proved)
``Deligne conjecture'' \cite{De4} which has fundamentally influenced modern
deformation theory.

Hence by Tamarkin's work $D_{\poly}^\Lscr(\Oscr)$ is
then also a \emph{strong homotopy Gerstenhaber algebra} ($G_\infty$-algebra
for short, see \S\ref{ref-4.3-8} below). Our main result will be the
following.
\begin{theorem}\label{ref-1.1-0} 
\begin{enumerate}
\item There exists a sheaf of $G_\infty$-algebras
  $\mathfrak{l}^\Lscr$ on $\Cscr$ together with
  $G_\infty$-quasi-isomorphisms
\begin{equation}\label{ref-1.1-1}
T_{\poly}^\Lscr(\Oscr)\longrightarrow\mathfrak{l}^\Lscr\longleftarrow D_{\poly}^\Lscr(\Oscr)\,.
\end{equation}
\item The isomorphism 
\[
T_{\poly}^\Lscr(\Oscr)\longrightarrow H^\ast(D_{\poly}^\Lscr(\Oscr))
\]
induced by \eqref{ref-1.1-1} is the usual
Hochschild-Kostant-Rosenberg isomorphism \cite{cal,ye2}. 
\end{enumerate}
\end{theorem}
The $G_\infty$-structure on  $D_{\poly}^\Lscr(\Oscr)$ depends on the
one time choice of a Drinfeld associator. Likewise the quasi-isomorphisms
in \eqref{ref-1.1-1} depend on the one time choice of a local formality
isomorphism. Once these choices are made the quasi-isomorphisms are canonical.

\medskip

In Appendix \ref{ref-A-39} we prove a strengthening of Halbout's globalization
result \cite[Theorem 4.5]{halb} for Tamarkin's local quasi-isomorphism. 

\medskip

The second author thanks Dmitry Tamarkin for some useful discussions. 

\section{Notations and conventions}
Our grading conventions for Gerstenhaber, $B_\infty$- and 
  $G_\infty$-algebras are shifted with respect to the usual ones. In our setup
the Lie bracket has degree zero and the cupproduct has degree one.

\section{Preliminaries on $B_\infty$-algebras}

\subsection{$B_\infty$-algebras}
\label{ref-3.1-2}

For a detailed discussion of $B_\infty$-algebras we refer to
\cite[\S5.2]{GJ}. Let $V$ be a graded vector space and $T^c(V)$ be the cofree 
cotensor algebra (with counit).  As graded vector spaces we have $T^c(V)=T(V)$.
The comultiplication is given by
\[
\Delta(a_1|\cdots|a_n)=\sum_i (a_1|\cdots |a_i)\otimes (a_{i+1}|\cdots|a_n)
\]
where
\[
(a_1|\cdots|a_n)\overset{\text{def}}{=}a_1\otimes\cdots\otimes a_n
\in T^{c,n}(V)
\]
and $1=()$. The counit is given by
\[
\epsilon(a_1|\cdots|a_n)=
\begin{cases}
1&\text{if $n=0$}\\
0&\text{otherwise}
\end{cases}
\]
A $B_\infty$-structure on $V$ consists of a DG-bialgebra
structure 
\[
(T^c(V),\Delta,\epsilon,m,1,Q)
\]
on $T^c(V)$ with unit equal to $1\in k=T^{c,0}(V)$. We also say that $V$
is $B_\infty$-algebra. A $B_\infty$-algebra morphism $V\rightarrow W$ is a morphism
of DG-bialgebras $\psi:T^c(V)\rightarrow T^c (W)$. If $\psi$ is obtained by extending
a morphism of DG-vector spaces $\psi:V\rightarrow W$ then we say that $\psi$ is strict. 

If $(T^c(V),\Delta,\epsilon,m,1,Q)$ is a $B_\infty$-structure on $V$ then
$Q$ is determined by the compositions
\[
Q^i:T^{c,i}(V)\hookrightarrow T^{c}(V)\xrightarrow{Q}
T^c(V)\xrightarrow{\text{projection}} V
\]
and likewise $m$ is determined by
\[
m^{p,q}:T^{c,p}(V)\otimes T^{c,q}(V)\hookrightarrow T^{c}(V)\otimes T^{c}(V)
\xrightarrow{m} T^c( V)\xrightarrow{\text{projection}} V
\]
A $B_\infty$-structure leads in a natural way to one unary and two
binary operations on~$V$ given respectively by a differential $Q^1$, a
Lie bracket $[v,w]=m^{1,1}(v,w)-(-1)^{|v||w|}m^{1,1}(w,v)$ of degree
zero and a ``cupproduct'' $Q^2$ of degree one. The computations in
\cite[\S 5.2]{GJ} show that $(V,Q^2,[-,-],Q^1)$ is a DG-Gerstenhaber
algebra up to homotopies expressible in the ternary operations $Q^3$,
$m^{1,2}$, $m^{2,1}$. In particular $H^\ast(V)$ is a Gerstenhaber
algebra.

\subsection{Inner $B_\infty$-algebras and brace algebras}

The $B_\infty$-algebras which we encounter below are of a
special kind.  Assume that $(T^c(V),\Delta,\epsilon,m,1)$ is a bialgebra. Let
$\mu\in V_1$ be such that $m(\mu,\mu)=0$. I.e. $m^{1,1}(\mu,\mu)=0$.
Since $\Delta(\mu)=\mu\otimes 1+1\otimes \mu$ we know that $Q=[\mu,-]$
is a degree one biderivation\footnote{In this context, a {\em biderivation} is a graded 
linear map that is both a derivation (for the product) and a coderivation 
(for the coproduct). } 
on $(T^c(V),\Delta,\epsilon,m,1)$. Hence this gives a
$B_\infty$-structure on $V$. We will say that $V$ is an \emph{inner}
$B_\infty$-algebra.

\medskip

Another major simplification appears when  $m^{p,q}=0$ for $p>1$. 
A $B_\infty$-algebra satisfying this condition is called a \emph{brace algebra}. 

It the case of a brace algebra it is easier to express the associativity
condition for $m$.  We write
\[
m^{1,p}(a,(b_1|\cdots |b_p))=a\{b_1,\cdots, b_p\}
\]
and the full multiplication can be expressed by the following identity
\[
m((a_1|\cdots|a_p),(b_1|\cdots|b_q))
=\sum_{0\le i_1\le \cdots\le i_p\le q}
(-1)^\epsilon (b_1|\cdots |b_{i_1}| a_1\{b_{i_1+1},\cdots\}|
\cdots |b_{i_p}| a_p\{b_{i_p+1},\cdots\}|\cdots |b_q)
\]
where $(-1)^\epsilon$ is the sign obtained from passing the $a$'s across the $b$'s.

The associativity condition becomes
\[
a\{b_1,\ldots,b_q\}\{c_1,\ldots,c_r\}
=
\sum_{0\le i_1\le \cdots\le i_q\le r}
(-1)^\epsilon a\{c_1,\ldots,c_{i_1}, b_1\{c_{i_1+1},\cdots\},
\cdots ,c_{i_q}, b_q\{c_{i_q+1},\cdots\},\cdots ,c_r\}
\]
where $\epsilon$ is a usual. 

\section{Preliminaries on $G_\infty$-algebras}
\label{ref-4-3}

\subsection{Co-strict homotopy Lie bialgebras}
\label{ref-4.1-4}
Assume that $(\frak{g},\delta)$ is a Lie coalgebra.
We view the coproduct
\[
\delta:\frak{g}\rightarrow \frak{g}\otimes \frak{g}
\]
as a map of degree $-1$
\[
\bar{\delta}:\frak{g}[1]\rightarrow \frak{g}[1]\otimes \frak{g}[1]
\]
via
\[
\bar{\delta}(g)=(-1)^{|g_{[1]}|}g_{[1]}\otimes g_{[2]}
\]
where $\delta(g)=g_{[1]}\otimes g_{[2]}$.
We can extend $\bar{\delta}$ to  a bicoderivation of degree -1
\[
\bar{\delta}:S(\frak{g}[1])\rightarrow S(\frak{g}[1])\otimes S(\frak{g}[1])
\]
via the formula 
\begin{equation}
\label{ref-4.1-5}
\bar{\delta}(g_1\cdots g_n)=\sum_{i}(-1)^{|g_1|+\cdots+|g_{i-1}|+i-1}
\Delta(g_1)\cdots \Delta(g_{i-1})\bar{\delta}(g_i)\Delta(g_{i+1})
\cdots \Delta(g_n)
\end{equation}
where $\Delta$ is the coshuffle coproduct on $S(\frak{g}[1])$ and
where $|g_j|$ refers to the degree of $g_i$ as an element of $\frak{g}$
(and not of $\frak{g}[1]$).  In this way $S(\frak{g}[1])$ becomes a
Gerstenhaber coalgebra.

A \emph{co-strict homotopy Lie bialgebra} (CSHLB for short) structure
on $\frak{g}$ is a coderivation $Q$ on $S(\frak{g}[1])$ of degree one
and square zero such that
\[
\bar{\delta}\circ Q+(Q\otimes{\rm
    id}+{\rm id}\otimes Q)\circ\overline{\delta}=0
\]
and $Q(1)=0$. Thus a CSHLB-structure is an $L_\infty$-structure which
satisfies a suitable compatibility with respect to the cobracket. 
A {\em CSHLB-morphism} $\mathfrak{g}\to\mathfrak{h}$ is an $L_\infty$-morphism 
commuting with $\bar{\delta}$. 
\begin{example} \label{ref-4.1-6} If $\frak{g}$ is a DG-Lie bialgebra then
  its induced $L_\infty$-structure makes it into a CSHLB. A similar
  statement is true for DG-Lie bialgebra morphisms.
\end{example}

\subsection{Cofree Lie coalgebras}

Let $V$ be a graded vector space. The shuffle product $m_s$ makes $T^c(V)$
into a bialgebra. The \emph{cofree Lie coalgebra} cogenerated by $V$
is defined by
\begin{equation}
\label{ref-4.2-7}
L^c(V)=\ker\epsilon/m_s(\ker \epsilon,\ker\epsilon)
\end{equation}
with cobracket
\[
\delta=\Delta-\Delta^\circ
\]
It is convenient to denote the image of $v_1\otimes\cdots \otimes v_n$
in $L^{c,n}(V)$ by $\underline{v_1\otimes\cdots \otimes v_n}$. 
\begin{remarks} The formula \eqref{ref-4.2-7} is dual to the realization
  of the free Lie algebra $L(V)$ inside the the free algebra $T(V)$ as
  the primitive elements for the coproduct given by
  $\Delta(v)=v\otimes 1+1\otimes v$ (which is dual to the shuffle
  product).
\end{remarks}

\subsection{$G_\infty$-algebras}
\label{ref-4.3-8}

A {\em $G_\infty$-structure} on $V$ is a CSHLB-structure on $L^c(V)$.
Likewise a $G_\infty$-morphism $V\rightarrow W$ is a CSHLB-morphism $L^c(V)\rightarrow
L^c(W)$. If this morphism is obtained by extending a morphism $V\rightarrow W$
of DG-vector spaces then we call it strict.
For a detailed study of $G_\infty$-algebras
see~\cite{gin}.

A $G_\infty$-structure is determined by maps 
\[
Q^{p_1,\ldots,p_n}: L^{c,p_1}(V)[1]\otimes
\cdots\otimes L^{c,p_n}(V)[1]\rightarrow V[1]
\]
of degree $1$ where $Q^n$ is the sum of all $Q^{p_1,\ldots,p_n}$. We obtain corresponding 
maps of degree $2-n$ 
\[
l^{p_1,\ldots,p_n}:L^{c,p_1}(V)\otimes\cdots\otimes L^{c,p_n}(V)\rightarrow V
\]
using standard sign conventions 
\[
l^{p_1,\ldots,p_n}(\underline{\gamma}_1,\ldots,\underline{\gamma}_n)=
(-1)^{n+(n-1)|\gamma_1|+(n-2)|\gamma_2|\cdots+|\gamma_{n-1}|}
Q^{p_1,\ldots,p_n}(\underline{\gamma}_1,\ldots,\underline{\gamma}_n)\,.
\]
Likewise a $G_\infty$-morphism $\psi:V\rightarrow W$ is determined by maps 
\[
\psi^{p_1,\ldots,p_n}: L^{c,p_1}(V)[1]\otimes
\cdots\otimes L^{c,p_n}(V)[1]\rightarrow W[1]
\]
of degree $0$. 

Recall that an $L_\infty$-structure on $V$ is given by a coderivation
of degree one and square zero on $S(V[1])$. It is determined by maps
\[
Q^n:S^n(V[1])\rightarrow V[1]
\]
of degree one, or likewise by maps
\[
m^n:\bigwedge^l V\rightarrow V
\]
of degree $2-l$. 
A $G_\infty$-algebra becomes an $L_\infty$-algebras by putting
\[
m^n=l^{1,1,\ldots,1}
\]
A $G_\infty$-structure leads in a natural way to one unary and two
binary operations on~$V$ given respectively by a differential $l^1$, a
Lie bracket up to homotopy $l^{1,1}$ of degree
zero and a ``cupproduct'' $l^2$ of degree one.

A computation shows that $(V,l^{1,1},l^2,l^1)$ satisfies the
Gerstenhaber axioms up to homotopies in terms of the ternary operations
$l^{3}$, $l^{2,1}$, $l^{1,1,1}$.  In particular
$H^\ast(V)$ is a Gerstenhaber algebra.

\section{Operads and Tamarkin's morphism}
\label{ref-5-9}

It is easy to see that $B_\infty$-algebras and $G_\infty$-algebras are
DG-algebras over certain DG-operads which we denote respectively by
$B_\infty$ and $G_\infty$.  The operad $B_\infty$ is generated as a graded
operad by the
operations $m^{p,q}\in B_\infty(p+q)$ (of degree zero)
and the operations $Q^n\in B_\infty(n)$, $n\ge 2$ (of degree one).

The operad $G_\infty$ is freely generated as a graded operad by the
operations $l^{p_1,\ldots,p_n}\in G_\infty(p_1+\cdots+p_n)$,
$(p_1,\ldots,p_n)\neq (1)$ (of degree $2-n$) and its homology is the
graded operad of Gerstenhaber algebras, which we denote by $G$.
Likewise in \S\ref{ref-3.1-2} we have indicated that a $B_\infty$-algebra
is a Gerstenhaber algebra up to homotopy.  This corresponds the
existence of a morphism of operads $G\rightarrow H^\ast(B_\infty)$ preserving
Lie bracket and cup product.

Tamarkin's
amazing discovery is the following.
\begin{theorem}[\cite{Tamarkin}]
There exists a morphism of DG-operads $T:G_\infty\rightarrow B_\infty$ such that the
following diagram is commutative 
\[
\begin{CD}
H^\ast(G_\infty) @>H^\ast(T)>> H^\ast(B_\infty)\\
@A\cong AA @AAA\\
G@= G
\end{CD}
\]
\end{theorem}
We discuss some further properties of $T$. For more details we refer to
\cite{Tamarkin}.  
\begin{lemma}
\label{ref-5.2-10}
One has $T(l^{p_1,\ldots,p_n})=0$ for $n>2$. Furthermore $T(l^{p_1,p_2})$
can be expressed in terms of operations $m^{p'_1,p'_2}$.
\end{lemma}
\begin{proof} This follows immediately by degree reasons. Indeed
  $|l^{p_1,p_2,\ldots,p_n}|<0$ and $B_\infty$ contains no elements of
  strictly negative degree. Similarly $T(l^{p_1,p_2})$ has degree zero
  and hence it must be expressible in terms of the generators
  $m^{p_1',p'_2}$ since the $Q^n$ have strictly positive degree.
\end{proof}
We have mentioned in \S\ref{ref-4.3-8} that any $G_\infty$-algebra is an $L_\infty$-algebra. This
corresponds to the morphism of operads $L_\infty\rightarrow G_\infty$ which
sends $m^n$ to $l^{1,\ldots,1}$ with the $1$'s occurring $n$ times. 

Let $L$ be the operad of Lie algebras. Then we also have the usual
quasi-isomorphism of operads $L_\infty\rightarrow L$ which kills $l^n$ for
$n>2$.

Finally we have a morphism of DG-operads $L\rightarrow B_\infty$ which sends $m^2$
to $m^{1,1}-\sigma(m^{1,1})$ where $\sigma=(12)$. 
\begin{lemma}[\cite{Tamarkin}] There is a commutative diagram of DG-operads
\[
\begin{CD}
G_\infty @>T>> B_\infty\\
@AAA @AAA\\
L_\infty @>>> L
\end{CD}
\]
\end{lemma}
\begin{proof} The proof is easy so for the benefit of the reader we
  will give it here. Since $T(l^{1,1,\ldots,1})=0$ the diagram is
  commutative when evaluated on $l^{1,1,\ldots,1}$.

  It remains to prove $T(l^{1,1})=m^{1,1}-\sigma(m^{1,1})$. Now
  $T(l^{1,1})$ is an element of $B_\infty(2)$ which has degree zero and
  is anti-symmetric. Since $B_\infty(2)_0=km^{1,1}+k\sigma(m^{1,1})$ we
  deduce $T(l^{1,1})=\alpha(m^{1,1}-\sigma(m^{1,1}))$ for $\alpha\in k$.
  Descending to cohomology yields $\alpha=1$.
\end{proof}
Let $\overline{G}_\infty$ be the quotient of $G_\infty$
by $l^{p_1,\ldots,p_n}$, $n>2$. By lemma \ref{ref-5.2-10} we deduce that $T$
factors as follows
\[
\overline{T}:\overline{G}_\infty\rightarrow B_\infty
\]
An $\overline{G}_\infty$-structure on $V$ corresponds to a DG-Lie bialgebra
structure
$(L^c(V),\delta,[-,-],Q)$ on $L^c(V)$ where $\delta$ is the standard
cobracket (cf.~Example \ref{ref-4.1-6}). So if $V$ is a $B_\infty$-algebra (e.g.\
a Hochschild complex) then the $G_\infty$-structure on $T_\ast(V)$ is actually 
quite special. 
\begin{remark} The theory of (de)quantization for Lie bialgebras
  \cite{EK1,EK2} shows that the operads $\overline{G}_\infty$ and
  $B_\infty$ are closely related \cite{GH}. However they are not
  isomorphic as
\begin{align*}
\overline{G}_\infty(2)_0&=kl^{1,1}\\
{B}_\infty(2)_0&=km^{1,1}\oplus k\sigma(m^{1,1})
\end{align*}
Thus $\overline{G}_\infty(2)_0$ is one dimensional and $B_\infty(2)_0$ is two
dimensional.
\end{remark}
We will need the following technical result which is a straightforward generalization of \cite[Lemma 4.7]{halb}. 
\begin{proposition}\label{ref-5.5-11}
Assume that $V$ is a $\overline{G}_\infty$-algebra. In particular $(V,l^{1,1})$ is a graded Lie
algebra. Assume that $U$ is a perfect Lie subalgebra of $V$ (i.e.\ $l^{1,1}(\wedge^2 U)=U$) such 
that for every $u\in U$, $l^{1,1}(u,-)$ acts as a derivation with respect to the operations $l^{1,q}$. 
Then $l^{1,q}(u,-)=0$ for all $u\in U$ and $q>1$.
\end{proposition}
\begin{proof}
  In the computations below the $u$'s are elements of $U$ and the
  $v$'s are elements of $V$. Let $(L^c(V),\delta,[-,-])$ be the Lie
bialgebra structure on $L^c(V)$ corresponding to the $\overline{G}_\infty$-structure. 

Let $l^{1,1}(u,-)$ act on $L^{c}(V)$ by
  extending the action on $V$ using the obvious Leibniz rule:
\[
l^{1,1}(u,\underline{v_1\otimes \cdots \otimes  v_q})=
\sum_i (-1)^{(|v_1|+\cdots+|v_{i-1}|)|u|} \underline{v_1\otimes\cdots
\otimes l^{1,1}(u,v_i)\otimes\cdots \otimes v_q}
\]  
Assume
  that we have shown that $l^{1,i}(u,-)=0$ for $1<i<q$ for $q\ge 2$ (for
$q=2$ there is nothing to show). We
  will prove that $l^{1,q}(u,-)=0$. A trite computation using the
induction hypotheses shows for $\gamma\in L^{q,c}(V)$
\[
[u,\gamma]=l^{1,1}(u,\gamma)+l^{1,q}(u,\gamma).
\]
Therefore,
  writing out the Jacobi identity
  $[u_1,[u_2,\gamma]]-(-1)^{|u_1||u_2|}[u_2,[u_1,\gamma]]-[[u_1,u_2],\gamma]=0$ yields
\begin{align*}
& l^{1,1}(u_1,l^{1,1}(u_2,\gamma))+l^{1,q}(u_1,l^{1,1}(u_2,\gamma))
+l^{1,1}(u_1,l^{1,q}(u_2,\gamma))
\\
& -(-1)^{|u_1||u_2|}l^{1,1}(u_2,l^{1,1}(u_1,\gamma))
-(-1)^{|u_1||u_2|}l^{1,q}(u_2,l^{1,1}(u_1,\gamma))
\\
& -(-1)^{|u_1||u_2|}l^{1,1}(u_2,l^{1,q}(u_1,\gamma))
-l^{1,1}(l^{1,1}(u_1,u_2),\gamma)-l^{1,q}(l^{1,1}(u_1,u_2),\gamma)=0\,.
\end{align*}
After projecting onto $V$ we obtain 
\begin{multline}\label{ref-5.1-12}
l^{1,q}(u_1,l^{1,1}(u_2,\gamma))
+l^{1,1}(u_1,l^{1,q}(u_2,\gamma))
-(-1)^{|u_1||u_2|}l^{1,q}(u_2,l^{1,1}(u_1,\gamma))\\
-(-1)^{|u_1||u_2|}l^{1,1}(u_2,l^{1,q}(u_1,\gamma))
)-l^{1,q}(l^{1,1}(u_1,u_2),\gamma)=0\,.
\end{multline}
According to our hypotheses $l^{1,1}(u,-)$ is a derivation with respect to 
the operation $l^{1,q}$. Thus
\[
l^{1,1}(u_1,l^{1,q}(u_2,\gamma))=l^{1,q}(l^{1,1}(u_1,u_2),\gamma)+(-1)^{|u_1||u_2|}l^{1,q}(u_2,l^{1,1}(u_1,\gamma))
\]
and thus we get $l^{1,q}(u_1,l^{1,1}(u_2,\gamma))=(-1)^{|u_1||u_2|}l^{1,1}(u_2,l^{1,q}(u_1,\gamma))$. 
Exchanging $u_1$ and $u_2$ we also have
$l^{1,q}(u_2,l^{1,1}(u_1,\gamma))=(-1)^{|u_1||u_2|}l^{1,1}(u_1,l^{1,q}(u_2,\gamma))$. 
Then substituting back in \eqref{ref-5.1-12} we get 
\[
l^{1,q}(l^{1,1}(u_1,u_2),\gamma)=0\,.
\]
Finally, we use the fact that $U$ is perfect to obtain $l^{1,q}(U,-)=0$. 
\end{proof}

\section{Twisting}

\subsection{Twisting for $B_\infty$-algebras}

Assume that $V$ is a $B_\infty$-algebra. I.e.\ we have a DG-bialgebra
$(T^c(V),\Delta,\epsilon,m,1,Q)$ with the usual coalgebra structure. Assume
that $\omega\in V_1$ is a solution of the following Maurer-Cartan equation
\[
Q(\omega)+m(\omega,\omega)=Q^1(\omega)+m^{1,1}(\omega,\omega)=0
\]
then 
\[
(T^c(V),\Delta,\epsilon,m,1,Q_\omega)
\]
with
$Q_\omega(\gamma)=Q+m(\omega,\gamma)-(-1)^{|\gamma|}m(\gamma,\omega)$ defines a new DG-bialgebra. 
We denote the resulting 
$B_\infty$-structure on $V$ by $V_\omega$. Thus explicitly the new $B_\infty$-algebra structure 
is given by
\begin{equation}\label{ref-6.1-13}
\begin{split}
m_\omega^{p,q}&=m^{p,q}\\
Q^p_\omega(\gamma)&=Q^p(\gamma)+m^{1,p}(\omega,\gamma)-(-1)^{|\gamma|}m^{p,1}(\gamma,\omega)
\end{split}
\end{equation}

\subsection{Twisting for $G_\infty$-algebras}
\label{ref-6.2-14}

Here we will work with infinite series. We assume that our vector
spaces are equipped with suitable topologies, that the series we use
converge, and furthermore that standard series manipulations are
allowed. All these hypotheses hold in our applications.
See \cite[\S4, \S6.2]{vdbcalaque} for a more precise
description of a setting in which these hypotheses hold.

Assume that $(\mathfrak{g},\delta,Q)$ is a CSHLB.
Let $\omega\in\mathfrak{g}_1$ be a solution of the
$L_\infty$-Maurer-Cartan equation in $\mathfrak{g}$, 
\[
\sum \frac{1}{i!}Q^n(\omega^i)=0
\]
and define $Q_\omega$ as in \cite{ye7}. I.e.
\begin{equation}
\label{ref-6.2-15}
Q_\omega^i(\gamma)=\sum_{j\ge 0} \frac{1}{j!} Q^{i+j}
(\omega^j \gamma)\qquad \text{(for $i>0$)}
\end{equation}
\begin{propositions}
\label{ref-6.2.1-16}
  If $\delta(\omega)=0$, then
  $\mathfrak{g}_\omega:=(\mathfrak{g},\delta, Q_\omega)$ is a CSHLB.
\end{propositions}
\begin{proof} 
It follows from \cite{ye7} that $(\frak{g},Q)$ is an
  $L_\infty$-algebra.  
So we only have to check that $\bar{\delta}\circ
  Q_\omega+(Q_\omega\otimes{\rm id}+{\rm id}\otimes
  Q_\omega)\circ\overline{\delta}=0$. 

We follow the method of \cite{ye7}. Let $S(\frak{g}[1])\,\hat{}$ be the
completion of $S(\frak{g}[1])$ at the ideal generated by $\frak{g}[1]$.
This is a topological Hopf algebra and the part cogenerated by $\frak{g}[1]$
is $S(\frak{g}[1])$. 
One has $e^\omega\in S(\frak{g}[1])\,\hat{}$ and furthermore \cite{ye7}
$e^{\omega}$ is group like. I.e.
\begin{equation}
\label{ref-6.3-17}
\Delta(e^{\omega})=e^{\omega}\otimes e^{\omega}
\end{equation}
In particular multiplication by $e^{\omega}$ defines a coalgebra automorphism
of $S(\frak{g}[1])\,\hat{}$. 

According to \cite{ye7} one has 
$Q_\omega(\gamma)=e^{-\omega}Q(e^{\omega}\gamma)$. 
It easily follows that one has to prove
\begin{equation}
\label{ref-6.4-18}
\bar{\delta}(e^{\omega}\gamma)=(e^\omega\otimes e^\omega )\bar{\delta}(\gamma)
\end{equation}
According to the explicit formula \eqref{ref-4.1-5} we have
\[
\bar{\delta}(\omega^n\gamma)=\Delta(\omega^n)\bar{\delta}(\gamma)
\]
(using $\delta(\omega)=0$) and hence
\[
\bar{\delta}(e^{\omega}\gamma)=\Delta(e^{\omega})\bar{\delta}(\gamma)
\]
Invoking \eqref{ref-6.3-17} finishes the proof. 
\end{proof}
Assume now that $\psi:(\mathfrak{g},\delta,Q)\to(\mathfrak{h},\delta',Q')$
is a CSHLB-morphism and define $\omega'$ and $\psi_\omega$ as in
\cite{ye7}.
\begin{align*}
\psi_\omega^i(\gamma)&=\sum_{j\ge 0} \frac{1}{j!} \psi^{i+j}
(\omega^j \gamma)\qquad \text{(for $i>0$)}\\
\omega'&=\sum_{j\ge 1} \frac{1}{j!} \psi^j
(\omega^j )
\end{align*}
By \cite{ye7} $\omega'$ is a solution to the $L_\infty$-Maurer Cartan
equation in $\frak{g}_\omega$.
\begin{propositions} Assume $\delta(\omega)=0$. Then one has
  $\delta'(\omega')=0$ and $\psi_\omega$ defines a 
  CSHLB-morphism from $\mathfrak{g}_\omega$ to
  $\mathfrak{h}_{\omega'}$.
\end{propositions}
\begin{proof} By \cite{ye7} one knows that $\psi_\omega$ is a
  $L_\infty$-morphism. We prove $\delta'(\omega')=0$.  From
  \eqref{ref-4.1-5} we deduce $\bar{\delta}(e^\omega)=0$. Applying
  $\psi$ and using $e^{\omega'}=\psi(e^{\omega})$ (see \cite{ye7}) we deduce
$\bar{\delta}'(e^{\omega'})=0$. Using \eqref{ref-4.1-5} again we find
$\bar{\delta}'(e^{\omega'})=\Delta(e^{\omega'})\bar{\delta}'(\omega')$. 
Since $\Delta(e^{\omega'})=e^{\omega'}\otimes e^{\omega'}$ is invertible
we conclude $\bar{\delta}'(\omega')=0$. 

To show that
  $\bar{\delta'}\circ\psi_\omega=(\psi_\omega\otimes\psi_\omega)\circ\bar{\delta}$
one uses $\psi_\omega(\gamma)=e^{-\omega}\psi(e^{\omega}\gamma)$ \cite{ye7}
together with \eqref{ref-6.4-18} (and the corresponding equation
for $\bar{\delta}'$).
\end{proof}
Assume now that $(S(L^c(V)[1]),Q)$ is a $G_\infty$-structure on $V$ and 
$\omega\in V_1$ is a solution to the $L_\infty$-Maurer-Cartan equation
\begin{equation}
\label{ref-6.5-19}
\sum \frac{1}{i!}Q^{\overbrace{1,\dots,1}^i}(\omega^i)=0
\end{equation}
Since the standard cobracket on $L^c(V)$ is zero on $V$ we obtain by Proposition
\ref{ref-6.2.1-16} a new
$G_\infty$-structure on $V$ given by $(S(L^c(V)[1]),Q_\omega)$. We denote this
new $G_\infty$-structure by $V_\omega$.  Using \eqref{ref-6.2-15} we deduce
\begin{equation}
\label{ref-6.6-20}
Q_\omega^{p_1\dots p_i}(\underline{\gamma}_1\cdots\underline{\gamma}_i)=
\sum_{j\ge 0}\frac{1}{j!} Q^{\overbrace{1,\dots,1}^j,p_1,\dots,p_i}
(\omega^j\underline{\gamma}_1\cdots\underline{\gamma}_i )\qquad\text{(for $i>0$)}
\end{equation}

Similarly if $\psi:V\rightarrow W$ is a $G_\infty$-morphism and $\omega\in V_1$ is a solution
to the $L_\infty$-Maurer-Cartan equation then we obtain a twisted $G_\infty$-morphism
$\psi_\omega:V_\omega\rightarrow W_{\omega'}$ where $\omega'$ and $\psi_\omega$ are given
by the formulas 
\begin{align}
\label{ref-6.7-21} \psi_{\omega}^{p_1,\dots,p_i}(\underline{\gamma}_1\cdots\underline{\gamma}_i)&=
\sum_{j\ge 0}\frac{1}{j!}\psi^{\overbrace{1,\dots,1}^j,p_1,\dots,p_i}
(\omega^j\underline{\gamma}_1\cdots\underline{\gamma}_i)\qquad \text{(for $i>0$)}\\
\label{ref-6.8-22} \omega'&=\sum_{j\ge 1}\frac{1}{j!}\psi^{\overbrace{1,\dots,1}^j}(\omega^j )
\end{align}

\section{Descent for $G_\infty$-algebras and morphisms}
\label{ref-7-23}

\subsection{Descent for algebras over DG-operads}
\label{ref-7.1-24}

Let $O$ be a DG-operad and let $\widetilde{O}$ be its underlying
graded operad.  If $V$ is an algebra over $O$ then an $S$-action on
$V$ is a family $(i_s)_{s\in S}$ of $\widetilde{O}$-derivations of degree
$-1$ of $V$. Put $L_s=di_s+i_sd$
and define
\[
V^S=\{v\in V\mid \forall s\in S:i_s(v)=L_s(v)=0\}
\]
A straightforward computation shows that $V$ is an algebra over $O$. 
\subsection{Descent for $G_\infty$-morphisms}
Let $V$, $W$ be two $G_\infty$-algebras and let $\psi:V\rightarrow W$ be a
$G_\infty$-morphism. 
Assume that we are given $S$-actions
$(i_s)_{s\in S}$ on $V$ and $W$. We say that $\psi$ commutes with these
actions if $i_s$ acts as a derivation with respect to $\psi^{p_1,\ldots,p_n}$
for all $s\in S$. This is clearly equivalent to
\begin{equation}
\label{ref-7.1-25}
[\bar{\imath}_s,\psi]=0
\end{equation}
where the $\bar{\imath}_s$ are the coderivations
commuting with $\bar{\delta}$ on $S(L^c(V)[1])$ and
$S(L^c(W)[1])$ obtained by extending $i_s$. 
\begin{lemmas} Assume that $(i_s)_{s\in S}$ commute with $\psi$. Then $\psi$
restricts to a $G_\infty$-morphism $V^S\rightarrow W^S$. 
\end{lemmas}
\begin{proof} It is sufficient to prove that $\bar{L}_s$ commutes with all
$\psi^{p_1,\ldots,p_n}$ or equivalently that
\begin{equation}
\label{ref-7.2-26}
[\bar{L}_s,\psi]=0
\end{equation}
where $\bar{L}_s$ is the extension of $L_s$ to a coderivation on $S(L^c(V)[1])$
and $S(L^c(W)[1])$ commuting with $\bar{\delta}$. Since $i_s$ commutes
with all $Q^{p_1,\ldots,p_n}$ (except perhaps with $Q^1$) we obtain
\begin{equation}
\label{ref-7.3-27}
\bar{L}_s=[Q,\bar{\imath}_s]
\end{equation}
Thus \eqref{ref-7.2-26} follows from
\eqref{ref-7.3-27} and \eqref{ref-7.1-25}.
\end{proof}

\subsection{Compatibility of descent with twisting of $G_\infty$-algebras}
\label{ref-7.3-28}

Here we make the same blanket hypotheses on series manipulations as in
\S\ref{ref-6.2-14}. Let $V$ be a $G_\infty$-algebra,
let $\omega\in V$ be a solution of the $L_\infty$-Maurer-Cartan equation
\eqref{ref-6.5-19} and define $Q_\omega$ as in \eqref{ref-6.6-20}. 
\begin{lemmas}
  Assume that $(i_s)_{s\in S}$ is an $S$-action on $V$.
  Assume in addition that for all $(p_1,\ldots,p_n)\neq ()$ we have
  $Q^{1,p_1,\ldots,p_n} ((i_v\omega) \underline{\gamma}_1\cdots \underline{\gamma}_n)=0$. 
  Then $(i_s)_{s\in S}$ is an $S$-action on $V_\omega$.
\end{lemmas}
\begin{proof}
  If we compute $i_s(Q^{p_1,\dots,p_n}_\omega(\underline{\gamma}_1\cdots \underline{\gamma}_n))$ then
  we see that $i_s$ behaves itself as a derivation with respect to
  $Q^{p_1,\dots,p_n}_\omega$ except for extra terms of the form
  $$
Q^{\overbrace{1,\dots,1}^j,{p_1,\dots,p_n}} (\omega\cdots (i_s\omega)\cdots \omega
\underline{\gamma}_1\cdots \underline{\gamma}_n)\,.
  $$
  These are zero by
  hypotheses.
\end{proof}

\subsection{Compatibility of descent with twisting of $G_\infty$-morphisms}

Let $\psi:V\rightarrow W$ be a morphism of $G_\infty$-algebras.  Let $\omega\in V$
be a solution of the $L_\infty$ Maurer-Cartan equation \eqref{ref-6.5-19}.
Define $\psi_\omega$, $\omega'$ as in \eqref{ref-6.7-21}\eqref{ref-6.8-22}. 
\begin{lemmas}  
\label{ref-7.4.1-29}
Let $(i_s)_{s\in S}$
act on $V$ and $W$ and assume that $\psi$  commutes with it.
Assume  in addition that for all $(p_1,\ldots,p_n)$ we have
$\psi^{1,p_1,\ldots,p_n} ((i_v\omega) \underline{\gamma}_1\cdots\underline{\gamma}_n)=0$. 
Then $(i_s)_{s\in S}$ commutes with the twisted $G_\infty$-map
$\psi_\omega:V_\omega\rightarrow W_{\omega'}$. 
\end{lemmas}
\begin{proof}
  If we compute $i_s(\psi^{p_1,\dots,p_n}_\omega(\underline{\gamma}_1\cdots \underline{\gamma}_n))$
  then we see that $i_s$ behaves itself as a derivation
  with respect to $\psi^{p_1,\dots,p_n}_\omega$ except for extra terms
  of the form 
$$
\psi^{\overbrace{1,\ldots,1}^j,{p_1,\dots,p_n} }(\omega\cdots (i_s\omega)\cdots\omega
\underline{\gamma}_1\cdots\underline{\gamma}_n)\,.
$$
  These are zero by hypotheses.
\end{proof}

\section{Poly-vector fields and poly-differential operators}
\label{ref-8-30}

We briefly recall some notations from \cite{vdbcalaque}. For more details
the reader is referred to loc.\ cit.

Let $(\Cscr,\Oscr)$ be a ringed site and let $\Lscr$ be a locally free
Lie algebroid of rank $d$. The enveloping algebra of
$\Lscr$ is denoted by $U\Lscr$.  Its right $\Oscr$-module structure is
\emph{defined} to be the same as its left structure.  As in
\cite{vdbcalaque} $T_{\poly}^{\Lscr}(\Oscr)$ is the Lie algebra of
$\Lscr$-poly-vector fields on $(\Cscr,\Oscr)$ \cite{cal}. I.e.\ it is
the graded vector space $(\wedge_{\Oscr} \Lscr)[1]$ equipped with the
graded Lie bracket obtained by extending the Lie bracket on $\Lscr$.
We equip $T_{\poly}^{\Lscr}(\Oscr)$ with the standard cupproduct
(which is of degree one with our shifted grading).  In this way
$T_{\poly}^{\Lscr}(\Oscr)$ becomes a sheaf of Gerstenhaber algebras.

Similarly $D_{\poly}^{\Lscr}(\Oscr)$ is the DG-Lie algebra of
$\Lscr$-poly-differential operators on $(\Cscr,\Oscr)$ \cite{cal}.
I.e.\ it is the graded sheaf $T_\Oscr(U\Lscr)[1]$ equipped with a
structure of an inner brace algebra
\begin{align*}
D\{E_1,\ldots,E_n\}=\sum_{i_1+\cdots+i_n=|P|-n+1}
(-1)^\epsilon 
({\rm id}^{\otimes i_1}\otimes\Delta^{|E_1|}\otimes\cdots\otimes\Delta^{|E_n|}\otimes{\rm id}^{\otimes i_n})(D)\cdot\\
\cdot(1^{\otimes i_1}\otimes E_1\otimes\cdots\otimes E_n\otimes1^{\otimes i_n})\,,
\end{align*}
where $\epsilon=\sum_s(|Q_s|-1)i_s$. The biderivation $Q$ on
$T_k(D_{\poly}^{\Lscr}(\Oscr))$ is given by $[\mu,-]$ where
$\mu=1\otimes 1\in D_{\poly}^{\Lscr}(\Oscr)_1$. One checks that
$Q^n=0$ for $n>2$. So the only non-vanishing $B_\infty$-operations are
the braces and the cupproduct $Q^2$.  The latter is the ordinary product on
$D_{\poly}^{\Lscr}(\Oscr)=T_{\Oscr}(U\Lscr)$ (up to a sign).

If $\Lscr$ is omitted from the notation we assume
$\Lscr=\Dscr\mathit{er}_1(\Oscr,\Oscr)$. In that case
$T_{\poly}(\Oscr)$ and $D_{\poly}(\Oscr)$ are the ordinary sheaves of
poly-vector fields and poly-differential operators. 

We also consider relative variants $T_{\poly,\Ascr}(\Bscr)$,
$D_{\poly,\Ascr}(\Bscr)$ of these notations where $\Ascr$, $\Bscr$ are
sheaves of commutative DG-algebras (and $\Bscr$ is a
DG-$\Ascr$-algebra). These will be self explanatory. In this case $Q$
is given by $d+[\mu,-]$ where $d$ is the differential on
$D_{\poly,\Ascr}(\Bscr)$ obtained from the DG-structure on
$\Ascr,\Bscr$.  It is still true that the braces and the cupproduct
are the only non-vanishing $B_\infty$-operations on
$D_{\poly,\Ascr}(\Bscr)$. Furthermore the differential on
$D_{\poly,\Ascr}(\Bscr)$ is of the form $d_{\text{tot}}=d+d_{\text{Hoch}}$ where
$d_{\text{Hoch}}=[\mu,-]$.
\section{Tamarkin's local formality morphism}
\label{ref-9-31}
Let $F=k[[t_1,\dots,t_d]]$. In \cite{Tamarkin} Tamarkin proved the existence
of $G_\infty$-quasi-isomorphism 
\[
\Psi:T_{\rm poly}(F)\to T_\ast(D_{\rm poly}(F))
\]
such that $\Psi^{1}$ is given by the HKR formula.

Moreover, one can construct this quasi-isomorphism in such a way that 
it has the following properties (see \cite[Theorem 4.5]{halb} or Theorem \ref{ref-A.1.1-40} below): 
\begin{itemize}
\item[(P4)] $\Psi^{1,\dots,1}(\underline{\gamma}_1\cdots\underline{\gamma}_n)=0$ for 
$\gamma_i\in T_{\poly}(F)_0$ and\footnote{Note that the case $n>2$
is automatic for degree reasons.} $n\ge 2$, 
\item[(P5)] $\Psi^{1,p_2,\dots,p_n}(\gamma\underline{\alpha}_2\cdots\underline{\alpha}_n)=0$ for $n\geq2$, 
$\gamma\in\mathfrak{gl}_d(k)\subset T_{\rm poly}(F)_0$ and $\underline{\alpha}_i\in L^{c,p_i}(T_{\rm poly}(F))$. 
\end{itemize}
\section{Proof of Theorem \ref{ref-1.1-0}}
Our proof parallels the proof of the corresponding $L_\infty$-result
in \cite{vdbcalaque}. We let $J\Lscr$ be the sheaf of
$\Lscr$-jetbundles on $(\Cscr,\Oscr)$. Thus
\[
J\Lscr=\invlim \uHom_{\Oscr}((U\Lscr)_{\le n},\Oscr)
\]
We show in \cite{vdbcalaque} that $J\Lscr$ is in a natural way a
$U\Lscr\otimes_k U\Lscr$ module.  Hence we have two commuting actions
of $\Oscr$ and $\Lscr$ on $J\Lscr$, depending on whether we embed them
in the first or second copy of $U\Lscr$. As in \cite{vdbcalaque} we
denote these two actions by $(\Oscr_1,\Lscr_1)$ and
$(\Oscr_2,\Lscr_2)$. We may view these as flat connections on $J\Lscr$
and in particular we obtain that $\wedge \Lscr_1^\ast\otimes_{\Oscr_1}
J\Lscr$ is in a natural way a $\wedge\Lscr_1^\ast$-DG-algebra (the
latter being a natural analogue of the De Rham complex).

As above let $F=k[[t_1,\ldots,t_d]]$. Then there exists a natural sheaf of
commutative DG-$\wedge\Lscr_1^\ast$-algebras $C^{\coord,\Lscr}$ (see
\cite[\S5.2]{vdbcalaque}) such that\footnote{From here on we generally work
with objects which are complete for some topology. Hence tensor products are
completed. See \cite[\S4.1]{vdbcalaque}
for a more precise discussion of our setting.}
\begin{equation}
\label{ref-10.1-32}
(C^{\coord,\Lscr}\otimes_{\wedge\Lscr_1^\ast}\wedge\Lscr_1^\ast\ctimes_{\Oscr_1}
 J\Lscr,d)=
(C^{\coord,\Lscr}\ctimes_{\Oscr_1} J\Lscr,d) \cong (C^{\coord,\Lscr}\ctimes_k F,d\otimes 1+\omega)
\end{equation}
where we have denoted the natural differentials by ``$d$'' and where $\omega$
is a solution of the Maurer-Cartan equation in the DG-Lie algebra 
\[
C^{\coord,\Lscr}\ctimes_k \Der_k(F,F)=
C^{\coord,\Lscr}\ctimes_k T_{\poly}(F)_0\subset C^{\coord,\Lscr}\ctimes_k D_{\poly}(F)_0
\]

\medskip

In \cite{vdbcalaque} we also considered a certain sub-DG
$\wedge\Lscr^\ast_1$-algebra $C^{\aff,\Lscr}$ of
$C^{\coord,\Lscr}$ which can for example be obtained by descent (see \S\ref{ref-7.1-24}). 
 More precisely for each $v\in
\frak{gl}_d(k)$ there exists a $\wedge\Lscr_1^\ast$-linear derivation
$i_v$ on $C^{\coord}$ (as a graded sheaf of algebras) of degree $-1$ such that
\begin{equation}
\label{ref-10.2-33}
C^{\aff}=(C^{\coord})^{\mathfrak{gl}_d(k)}
\end{equation}
We now construct some strict $B_\infty$-morphisms (see \cite{vdbcalaque})
\[
D_{\poly}^{\Lscr_2}(\Oscr_2)\xrightarrow[\text{qi}]{\alpha}
D_{\poly,C^{\aff,\Lscr}}(C^{\aff,\Lscr}\ctimes_{\Oscr_1} J\Lscr)
\xrightarrow{\beta}
D_{\poly,C^{\coord,\Lscr}}(C^{\coord,\Lscr}\ctimes_{\Oscr_1} J\Lscr)
\xrightarrow[\cong]{\gamma}
(C^{\coord,L} \ctimes_k D_{\poly}(F))_{\omega}
\]
The map $\alpha$ is constructed by letting $\Oscr_2,\Lscr_2$ act on
$J\Lscr$. This is possible since these actions commute with the
$\Oscr_1,\Lscr_1$-actions (which were the only ones we used so far).
It has been shown in \cite{vdbcalaque} that $\alpha$ is a quasi-isomorphism.

The map $\beta$ is obtained by extending scalars.  It remains to
discuss the map $\gamma$. Using \eqref{ref-10.1-32} we obtain an
isomorphism
\[
D_{\poly,C^{\coord,\Lscr}}(C^{\coord,\Lscr}\ctimes_{\Oscr_1} J\Lscr)
\cong D_{\poly,C^{\coord,\Lscr}}(C^{\coord,\Lscr}\ctimes F)
\]
which commutes with cupproduct and braces and hence it is an
isomorphism as sheaves of $\widetilde{B}_\infty$ algebras ($\widetilde{B}_\infty$
is the underlying graded operad for $B_\infty$, see \S\ref{ref-7.1-24}). The
Hochschild differential on the left is sent to the Hochschild
differential on the right. The natural differential $[d,-]$ on the
left is sent to $[d+\omega,-]$ on the right.

Thus we get as sheaves of $B_\infty$-algebras.
\begin{equation}
\label{ref-10.3-34}
D_{\poly,C^{\coord,\Lscr}}(C^{\coord,L}\ctimes_{\Oscr_1} JL)
\cong (D_{\poly,C^{\coord,\Lscr}}(C^{\coord,\Lscr}\ctimes F),
d_{\text{tot}}+[\omega,-])
\end{equation}
Since $\omega\in C^{\coord,L} \ctimes_k D_{\poly}(F)_0$ we have for $q>1$
\begin{equation}
\label{ref-10.4-35}
\begin{split}
m^{1,q}(\omega,-)&=0\\
m^{q,1}(-,\omega)&=0.
\end{split}
\end{equation} 
A simple computation using \eqref{ref-6.1-13} and \eqref{ref-10.4-35}  yields as $B_\infty$-algebras
\[
D_{\poly,C^{\coord,L}}(C^{\coord,L}\ctimes F)_\omega=(D_{\poly,C^{\coord,L}}(C^{\coord,L}\ctimes F),d_{\text{tot}}+[\omega,-])
\]
and combining this with \eqref{ref-10.3-34} we obtain the strict $B_\infty$-isomorphism
$\gamma$. 

We now apply the functor $T_\ast$. We get strict maps of $G_\infty$-algebras:
\begin{multline*}
T_\ast(D_{\poly}^{\Lscr_2}(\Oscr_2))\xrightarrow[\text{qi}]{\alpha}
T_\ast(D_{\poly,C^{\aff,\Lscr}}(C^{\aff,\Lscr}\ctimes_{\Oscr_1} J\Lscr))
\xrightarrow{\beta}
T_\ast(D_{\poly,C^{\coord,\Lscr}}(C^{\coord,\Lscr}\ctimes_{\Oscr_1} J\Lscr))\\
\xrightarrow[\cong]{\gamma}
T_\ast((C^{\coord,L} \ctimes D_{\poly}(F))_{\omega})=
T_\ast(C^{\coord,L} \ctimes D_{\poly}(F))_\omega\cong
(C^{\coord,L}\ctimes T_\ast(D_{\poly}(F))_\omega
\end{multline*}
The equality is an instance of the compatibility of $T_\ast$ with
twisting.  It seems this is in general a subtle issue on which we will
come back in a future paper. However in this special case we can use
an {\it ad hoc} argument.

\begin{lemma} We have
\begin{align*}
T_\ast((C^{\coord,L} \ctimes D_{\poly}(F))_{\omega})&=
T_\ast(C^{\coord,L} \ctimes D_{\poly}(F), d_{\text{tot}}+[\omega,-])\\
&=
(\widetilde{T}_\ast(C^{\coord,L} \ctimes D_{\poly}(F)),d_{\text{tot}}+[\omega,-])\\
&=(\widetilde{T}_\ast(C^{\coord,L} \ctimes D_{\poly}(F)),d_{\text{tot}})_\omega\\
&=T_\ast(C^{\coord,L} \ctimes D_{\poly}(F))_\omega
\end{align*}
\end{lemma}
\begin{proof} The first equality has already been established. The
second and fourth equalities are tautologies. Hence it remains
to establish the third equality. 

From \eqref{ref-10.4-35} we deduce that
$[\omega,-]=m^{1,1}(\omega,-)-m^{1,1}(-,\omega)$ behaves as a derivation
with respect to the operations $m^{p,q}$.
According to Lemma \ref{ref-5.2-10} the operations $l^{p_1,p_2}$ can be
expressed in terms of the $m^{p,q}$. Hence $[\omega,-]$ acts as a
derivation with respect to the operations $l^{p_1,p_2}$. It is easy to see
that the Lie algebra $D_{\poly}(F)_0$ is perfect and hence the same holds
for $C^{\coord,\Lscr}\ctimes_k D_{\poly}(F)_0$ which contains $\omega$.  
We deduce from
Proposition \ref{ref-5.5-11} that $l^{1,q}(\omega,-)=0$ for $q>1$.
By Lemma \ref{ref-5.2-10} we also know that $l^{p_1,\dots,p_n}=0$ for $n\ge 3$.
Substituting all this information in \eqref{ref-6.6-20} we deduce
\begin{align*}
  Q^{p_1}_\omega(\underline{\gamma})&=Q^{p_1}(\underline{\gamma})&
  (\text{if $p_1>1$})\\
  Q^{1}_\omega(\underline{\gamma})&=Q^1(\underline{\gamma})+Q^{1,1}(\omega,\underline{\gamma})\\
Q^{p_1,\dots,p_n}(\underline{\gamma_1}\cdots \underline{\gamma_i})&=0&
(\text{if $n\ge 2$})
\end{align*}
Translating this back in terms of the $l$'s we see that twisting by
$\omega$ does nothing except adding $l^{1,1}(\omega,-)=[\omega,-]$ to
the underlying differential. This is precisely the content of the
third equality in the statement of the lemma. 
\end{proof}
We will now construct similar morphisms of sheaves of DG-Gerstenhaber
algebras
\[
T_{\poly}^{\Lscr_2}(\Oscr_2)\xrightarrow[\text{qi}]{\alpha'}
T_{\poly,C^{\aff,\Lscr}}(C^{\aff,\Lscr}\ctimes_{\Oscr_1} J\Lscr)
\xrightarrow{\beta'}
T_{\poly,C^{\coord,\Lscr}}(C^{\coord,\Lscr}\ctimes_{\Oscr_1} J\Lscr)
\xrightarrow[\cong]{\gamma'}
(C^{\coord,L} \ctimes_k T_{\poly}(F))_{\omega}
\]
Again  the map $\alpha'$ is constructed by letting $\Oscr_2,\Lscr_2$ act on
$J\Lscr$ and the map $\beta'$ is obtained by extending scalars.  It remains to
discuss the map $\gamma'$.
Using \eqref{ref-10.1-32} we obtain an
isomorphism
\[
T_{\poly,C^{\coord,\Lscr}}(C^{\coord,\Lscr}\ctimes_{\Oscr_1} J\Lscr)
\cong T_{\poly,C^{\coord,\Lscr}}(C^{\coord,\Lscr}\ctimes F)
\]
which commutes with the cupproduct and the Lie bracket and hence is an
isomorphism as sheaves of $\widetilde{G}$ algebras.  The natural
differential $[d,-]$ on the left is sent to $[d+\omega,-]$ on the
right. Thus we get as sheaves of $G$-algebras
\begin{equation}
\label{ref-10.5-36}
T_{\poly,C^{\coord,\Lscr}}(C^{\coord,L}\ctimes_{\Oscr_1} JL)
\cong (T_{\poly,C^{\coord,\Lscr}}(C^{\coord,\Lscr}\ctimes F),d+[\omega,-])
\end{equation}
Since $T_{\poly,C^{\coord,\Lscr}}(C^{\coord,\Lscr}\ctimes F)$ is a
Gerstenhaber algebra the only operations that are non-zero are
$l^{1,1},l^2$. From the formula \eqref{ref-6.6-20} we deduce that the only
effect of twisting by $\omega$ is changing the differential into
$d+[\omega,-]$. Thus we obtain
\[
T_{\poly,C^{\coord,\Lscr}}(C^{\coord,\Lscr}\ctimes F)_\omega=
 (T_{\poly,C^{\coord,\Lscr}}(C^{\coord,\Lscr}\ctimes F),d+[\omega,-])
\]
Combining this with \eqref{ref-10.5-36} we obtain $\gamma'$.

We will  now construct a commutative diagram of $G_\infty$-morphism
\begin{equation}
\label{ref-10.6-37}
\hspace*{-2cm}\xymatrix{
  \scriptstyle T_{\poly}^{\Lscr_2}(\Oscr_2)\ar[r]^-{\alpha'} &
  \scriptstyle T_{\poly,C^{\aff,\Lscr}}
(C^{\aff,\Lscr}\ctimes_{\Oscr_1} J\Lscr)\ar[r]^-{\beta'}
  \ar[d]^{\Phi^{\aff}} &
  \scriptstyle T_{\poly,C^{\coord,\Lscr}}(C^{\coord,\Lscr}\ctimes_{\Oscr_1} J\Lscr) 
  \ar[r]^-{\gamma'}\ar[d]^{\Phi^{\coord}}&
  \scriptstyle (C^{\coord,\Lscr}\ctimes T_{\poly}(F))_{\omega}
  \ar[d]^{(\Id\otimes \widetilde{\Psi})_\omega}\\
  \scriptstyle T_\ast(D_{\poly}^{\Lscr_2}(\Oscr_2))\ar[r]_-\alpha &
  \scriptstyle T_\ast(D_{\poly,C^{\aff,\Lscr}}(C^{\aff,\Lscr}\ctimes_{\Oscr_1} J\Lscr))\ar[r]_-\beta &
  \scriptstyle 
  T_\ast(D_{\poly,C^{\coord,\Lscr}}(C^{\coord,\Lscr}\ctimes_{\Oscr_1} J\Lscr)) \ar[r]_-{\gamma}&
  \scriptstyle (C^{\coord,\Lscr} \ctimes T_\ast(D_{\poly}(F)))_{\omega}
}
\end{equation}
in which the horizontal arrows are strict. $\Phi^{\coord}$ is defined by
\[
\Phi^{\coord}=\gamma^{-1}\circ (\Id\otimes \widetilde{\Psi})_\omega\circ
\gamma^{\prime}
\] 
and $\Phi^{\aff}$ is derived from $\Phi^{\coord}$ using descent.  
As graded sheaves we have
\begin{align}
  T_\ast(D_{\poly,C^{\aff,\Lscr}} (C^{\aff,\Lscr}\ctimes_{\Oscr_1} J\Lscr))
&=C^{\aff,\Lscr} \ctimes_{\Oscr_1} D_{\poly,\Oscr_1}(J\Lscr)\\
\label{ref-10.8-38}  T_\ast(D_{\poly,C^{\coord,\Lscr}} (C^{\coord,\Lscr}\ctimes_{\Oscr_1} J\Lscr))
&=C^{\coord,\Lscr} \ctimes_{\Oscr_1} D_{\poly,\Oscr_1}(J\Lscr)
\end{align}
As indicated above there is a $\frak{gl}_d(k)$ action
$C^{\coord,\Lscr}$.  If $v\in \frak{gl}_d(k)$ acts by $i_v$ then we
obtain a $ \frak{gl}_d(k)$-action on $D_{\poly,C^{\coord,\Lscr}}
(C^{\coord,\Lscr}\ctimes_{\Oscr_1} J\Lscr)$ as sheaf of
$B_\infty$-algebras by letting $v$ act as $[i_v,-]$. Under the
isomorphism \eqref{ref-10.8-38} this action can also be viewed as the
linear extension of the action on $C^{\coord,\Lscr}$.  Since the
operations of $\widetilde{G}_\infty$ can be expressed in those of
$\widetilde{B}_\infty$ via $\tilde{T}$ it is clear that we obtain a $
\frak{gl}_d(k)$-action on $ T_\ast(D_{\poly,C^{\aff,\Lscr}}
(C^{\aff,\Lscr}\ctimes_{\Oscr_1} J\Lscr))$ as $G_\infty$-algebra.  So
the formalism of \S\ref{ref-7-23} applies.

Combining \eqref{ref-10.8-38} with \eqref{ref-10.2-33} we obtain
\[
D_{\poly,C^{\aff,\Lscr}} (C^{\aff,\Lscr}\ctimes_{\Oscr_1} J\Lscr)=
D_{\poly,C^{\coord,\Lscr}} (C^{\coord,\Lscr}\ctimes_{\Oscr_1} J\Lscr)^{\mathfrak{gl}_d(k)}
\]
For similar reasons the formalism of \S\ref{ref-7-23} applies
also to $T_{\poly,C^{\aff,\Lscr}} (C^{\aff,\Lscr}\ctimes_{\Oscr_1}
J\Lscr)$. We omit the trivial details.

From (P4) we deduce that $\omega'=\omega$. It remains to show that
$\Phi^{\coord}$ commutes with the $\frak{gl}_d(k)$-action.  Then we
may as well show that $(\Id\otimes\Psi)_\omega$ commutes with the
$\frak{gl}_d(k)$-action. To this end we can use the criterion given by
Lemma \ref{ref-7.4.1-29}. We need to prove that $(\Id\otimes
\Psi^{1p_1,\dots,p_n})(i_v\omega \underline{\gamma}_1\cdots\underline{\gamma}_n)=0$. 
By \cite[Lemma 5.2.1]{vdbcalaque} we have
$i_v\omega=1\otimes v$.  It now suffices to invoke (P5). 

As in \cite{vdbcalaque} we put
$\mathfrak{l}^{\Lscr}=T_\ast(D_{\poly,C^{\aff,\Lscr}}
(C^{\aff,\Lscr}\ctimes_{\Oscr_1} J\Lscr))$.  Now note that if we
restrict to the $L_\infty$-setting then the maps we have constructed
here are the same as those in \cite{vdbcalaque}. Hence the fact
that we obtain quasi-isomorphisms 
\[
T_{\poly}^\Lscr(\Oscr)\xrightarrow{\Phi^{\aff}\circ\alpha'}\mathfrak{l}^\Lscr
\xleftarrow{\alpha} T_\ast(D_{\poly}^\Lscr(\Oscr))\,.
\]
as well as part (2) of Theorem \ref{ref-1.1-0} follows from \cite[Thm 7.4.1]{vdbcalaque}.
\appendix
\section{Proof of the globalization properties}
\label{ref-A-39}
\subsection{Introduction and statement of the main results}
Since properties (P4)(P5) (see \S\ref{ref-9-31}) are crucial
for us, and since Halbout's proof of them is somewhat sketchy (see
\cite[Theorem 4.5]{halb}) we give the proof below. This appendix
can be read more or less independently from the main paper. 

\medskip

Our main result is a generalization of Halbout's result. 
\begin{theorems} 
\label{ref-A.1.1-40}
Let $W$ be a vector space of dimension $d$. Put $R=SW$ (thus
$\Spec R=W^\ast$).  Let $\Aff(W)$ be the affine group associated to
$W^\ast$ and let $\afff(W)$ be its Lie algebra. We let $\Aff(W)$ act
on $R$ in the natural way.

Equip $D_{\poly}(R)$ with its functorial $G_\infty$-structure obtained
from a morphism of DG-operads $G_\infty\rightarrow B_\infty$ (see
\S\ref{ref-5-9}). Then there exists a $G_\infty$-quasi-isomorphism
$\Psi:T_{\poly}(R)\rightarrow D_{\poly}(R)$ with the following properties
\begin{enumerate}
\item[(P1)]  (The ``locality'' property.) The 
  $\Psi^{p_1,\ldots,p_n}:T_{\poly}(R)^{\otimes \sum_ip_i}\rightarrow
  D_{\poly}(R)$ are $R$-poly-differential operators (where we
  consider $T_{\\poly}(R)$ as a free $R$-module of finite rank).
\item[(P2)] $\Psi^1:T_{\\poly}(R)\rightarrow D_{\poly}(R)$ is the HKR
  quasi-isomorphism.  
\item[(P3)] $\Psi$ is equivariant for $\Aff(W)$.
\item[(P4)]
  $\Psi^{1,\ldots,1}(\underline{\gamma}_1,\ldots,\underline{\gamma}_n)=0$
  for $\gamma_i\in T_{\poly}(R)_0$ and $n\ge 2$.
\item[(P5)] $\Psi^{1,p_2,\dots,p_n}(\gamma\underline{\alpha}_2\cdots\underline{\alpha}_n)=0$ for $n\geq2$, 
$\gamma\in\afff(W)\subset T_{\rm poly}(R)_0$ and $\underline{\alpha}_i\in L^{c,p_i}(T_{\rm poly}(R))$. 
\end{enumerate}
\end{theorems}
Note that the conditions (P1)(P3) imply that $\Psi^{p_1\cdots p_n}$ can be
written as a sum over graphs just like Kontsevich's standard
$L_\infty$-quasi-isomorphism \cite{Ko3}. See \S\ref{ref-A.7-77} and
also \cite{GamHal,Shoikhet}.

\medskip

The proof of Theorem \ref{ref-A.1.1-40} will be given at the end of
\S\ref{ref-A.6-65}. It depends on a detailed analysis of the proof of
Tamarkin's theorem \cite{Tamarkin}. Recall that the latter asserts the
existence of a $G_\infty$-quasi-isomorphism $\Psi$ satisfying (P2).

\medskip

Here is a more detailed sketch of the proof. We follow  more
or less the same strategy as Halbout in \cite{halb}. Put
$\frak{g}=D_{\poly}(R)$ and $\frak{h}=T_{\poly}(R)$ and let $Q_0$, $Q_2$
be the codifferentials on $S^c(L^c(\frak{h})[1])$ and
$S^c(L^c(\frak{g})[1])$ representing the $G_\infty$-structures on
$\frak{h}$ and $\frak{g}$.

The theory of minimal models shows that there is a
$G_\infty$-structure $Q_1$ on $\frak{h}\overset{\text{HKR}}{\cong}H^\ast(\frak{g})$ and a
$G_\infty$-quasi-isomorphism
\[
\Psi':(\frak{h},Q_1)\rightarrow (\frak{g},Q_2)
\]
such that $\Psi^{\prime 1}=\text{HKR}$. 

The construction of such a minimal model can be made very explicit
using trees (see \S\ref{ref-A.2-43} and also
\cite{CL,KS2,Merkulov}). From this construction it follows easily
that $\Psi'$ satisfies (P1-5) and $Q_1$ is given by affine invariant
$R$-poly-differential operators.

Now $\frak{h}$ equipped with the usual zero differential is a
(shifted) Schouten algebra and Tamarkin shows that it is in fact
rigid. This means that there is a $G_\infty$-isomorphism
\[
\Psi'':(\frak{h},Q_0)\rightarrow (\frak{h},Q_1) 
\]
such that $\Psi^{\prime\prime 1}=\Id_{\frak{h}}$.  For the proof of
Tamarkin's theorem it now suffices to take $\Psi=\Psi'\Psi''$.

\medskip

The crucial rigidity result is proved by showing that the
corresponding deformation complex
\[
(\coder (S^c(L^c(\frak{h})[1]),S^c(L^c(\frak{h})[1]))^\delta,[Q_0,-])=
(\Hom(S^c(L^c(\frak{h})[1]),\frak{h}[1]),[Q_0-])
\]
is essentially acyclic (by the superscript $\delta$ we mean that we
only consider coderivations compatible with $\delta$, see \S\ref{ref-4.1-4}). In fact its
sole cohomology is a copy of $k$ located in degree~$-2$.

\medskip

Since we want $\Psi$ to satisfy stronger conditions we have
to choose $\Psi''$ more carefully. This amounts to working with smaller
deformation complexes. It turns out (P4) follows from (P1)(P3)(P5) so we will
not worry about it in this introduction. 

To take care of (P1)(P3)(P5) we have to use the complex
\begin{equation}
\label{ref-A.1-41}
 (\Diff(S^c(L^c(\frak{h})/\afff(W)[1]),\frak{h}[1]),[Q_0,-])^{\Aff(W)}
\end{equation}
where ``$\Diff$'' means we only consider cochains which yield
$R$-poly-differential operators $L^{c,p_1}(\frak{h})\otimes \cdots
\otimes L^{c,p_n}(\frak{h})\rightarrow \frak{h}$ for all $(p_i)_i$.  We prove the following
result.
\begin{theorems} 
\label{ref-A.1.2-42} We have 
\[
H^\ast(\Diff(S^c(L^c(\frak{h})/\afff(W)[1]),\frak{h}[1])^{\Aff(W)},[Q_0,-])\cong
S(\frak{gl}(W)[-2])^{\GL(W)}[2]
\]
\end{theorems}
Thus the complex \eqref{ref-A.1-41} is not acyclic. Luckily
we find that its cohomology lives only in even degree and as the
obstruction against rigidity lives in odd degree we are saved.

\medskip

Below we keep the notations introduced in this introduction. If $V$ is
a $\GL(W)$-representation then if confusion is possible we denote the
corresponding $\Aff(W)$-representation by $\overline{V}$ (thanks to the 
projection $\Aff(W)\to \GL(W)\cong\Aff(W)/\textrm{translation}$). E.g. as an 
$\Aff(W)$-module $T_{\rm poly}(R)=R\otimes_k\wedge\overline{W}^*$. 

\subsection{Minimal models and trees}
\label{ref-A.2-43}

Minimal models for strong homotopy algebras over operads can be
constructed using trees (see \cite{CL,KS2,Merkulov}). This approach
was popularized by Kontsevich and Soibelman for $A_\infty$-algebras in
\cite{KS2}. Below we give a straightforward account of this theory,
which shows for example that the side conditions sometimes imposed
on the homotopy are unnecessary.

Let $\Oscr$ be a graded operad.  Thus we have maps
\begin{equation}
\label{ref-A.2-44}
\Oscr(m)\otimes \Oscr(n_1)\otimes \cdots \otimes \Oscr(n_m)\r
\Oscr(n_1+\cdots +n_m)
\end{equation}
For simplicity we assume $\Oscr(0)=0$ (no
constants), $\Oscr(1)=k$ (no non-trivial unary operations). We also
assume $\dim \Oscr(n)<\infty$ for all $n$. 

A  coalgebra $C$ over $\Oscr$ is a collection of 
$S_n$-invariant ``operations''
\[
\Delta_m:\Oscr(m)\otimes C\rightarrow C^{\otimes m}
\]
satisfying the standard axioms.

The cofree $\Oscr$-coalgebra over the graded vector space $V$ is given
by the formula\footnote{For cofree coalgebras it would seem more natural to use a definition in terms of $S_n$-\emph{invariants}. However in characteristic zero
there is a canonical identification between invariants and coinvariants.}
\begin{equation}
\label{ref-A.3-45}
F^c(V)=\bigoplus_n\Oscr(n)^\ast \otimes_{S_n} V^{\otimes n}
\end{equation}
We remind the reader of the formula for the $\Oscr$-coaction. 

For $\phi\in \Oscr(n_1+\cdots+n_m)$ write
\[
\sum\phi_{(0)}\otimes\phi_{(1)}\otimes\cdots\otimes \phi_{(m)}
\]
for the image of $\phi$ under the map
\[
\Oscr(n_1+\cdots+n_m)^\ast\rightarrow \Oscr(m)^\ast\otimes \Oscr(n_1)^\ast
\otimes \cdots\otimes  \Oscr(n_m)^\ast
\]
dual to \eqref{ref-A.2-44}.
Then, for any $\gamma\in\mathcal O(m)$ and $\phi\in\mathcal O(n)^*$, we have 
\begin{multline*}
\Delta_m(\gamma\otimes \phi\otimes v_1\otimes\cdots \otimes v_{n})=
\sum_{\sigma\in S_n, \sum_{i=1}^m\!n_i=n}\pm
\frac{(\sigma\phi)_{(0)}(
\gamma)}
{n_1!\cdots n_m!} \bigl((\sigma\phi)_{(1)}\otimes v_{\sigma^{-1}(1)}
\otimes\cdots \otimes v_{\sigma^{-1}(n_1)}\bigr)
\otimes\\
\cdots
\otimes 
\bigl ((\sigma\phi)_{(m)}\otimes v_{\sigma^{-1}(n_1+\cdots+n_{m-1}+1)}
\otimes\cdots \otimes v_{\sigma^{-1}(n)}\bigr)
\end{multline*}
where here and below the sign is controlled by the Koszul
convention.\footnote{In such formulas it is possible to get rid of the inverse 
factorials by restricting the summation to shuffles.}

For use below we give the
formula for a coalgebra morphism $\psi:F^{c}(V)\rightarrow F^c(W)$ determined by
maps
\[
\psi^n:F^{c,n}(V)\rightarrow W
\]
\begin{multline}
\label{ref-A.4-46}
  \psi(\phi\otimes v_1\otimes\cdots\otimes v_n)
  =
  \sum_{\begin{smallmatrix}m\geq2,\sigma\in S_n\\ \sum_{i=1}^m n_i=n\end{smallmatrix}}\pm
  \frac{(\sigma\phi)_{(0)}}{m!\,n_1!\cdots n_m!}\otimes \psi^{n_1}\bigl((\sigma\phi)_{(1)}\otimes
  v_{\sigma^{-1}(1)}\otimes\cdots \otimes
  v_{\sigma^{-1}(n_1)}\bigr)\otimes\\\cdots \otimes
  \psi^{n_m}\bigl((\sigma\phi)_{(m)}\otimes
  v_{\sigma^{-1}(n_1+\cdots+n_{m-1}+1)}\otimes\cdots \otimes
  v_{\sigma^{-1}(n)}\bigr)
\end{multline}
Now we discuss minimal models.  Let $Q_2:F^c(V)\rightarrow V[1]$ be a
codifferential and let $H^\ast(V)$ be the cohomology of the complex
$(V,Q^1_2)$. Then $Q_2^1Q_2^2+Q_2^2Q_2^1=0$ and hence $H^\ast(Q_2^2)$
is well defined and yields a map $F^{c,2}(H^\ast(V))\rightarrow H^\ast(V)$.

Furthermore we have $Q^3_2Q^1_2+(Q^2_2)^2+Q^1_2Q^3_2=0$. Thus $(Q^2_2)^2$
is homotopic to the zero map. Hence $H^\ast(Q_2^2)^2=0$.

For a series of maps $\psi^i:F^{c,i}(H^\ast(V))\rightarrow V$, $i=1,\ldots,n$
we let $\psi^{\le n}$ be the corresponding coalgebra map
$F^c(H^\ast(V))\rightarrow F^c(V)$ (the higher Taylor coefficients are assumed
to be zero). We use similar conventions for coderivations. 

We now choose a decomposition $(V,Q^1_2)=(H^\ast(V),0)\oplus W$ as complexes
and we let $i:H^\ast(V)\rightarrow V$, $p:V\rightarrow H^\ast(V)$ be the corresponding
injection and projection map. Furthermore we choose a homotopy
$H:V\rightarrow V[-1]$ such that 
\[
ip-\Id=Q^1_2H+HQ^1_2
\]
Let $P:F^{c}(V)\rightarrow V$ be the projection and let $I_n:F^{n,c}(V)\r
F^c(V)$ be the injection. We use the same notation with $V$ replaced
by $H^\ast(V)$. 

\begin{propositions}
\label{ref-A.2.1-47}
 Define a coderivation $Q_1$ of degree one on
$F^{c,n}(H^\ast(V))$ and a graded  coalgebra
map $\psi:F^c(H^\ast(V))\rightarrow F^c(V)$ recursively as follows: $\psi^1=i$,
$Q^1_1=0$ and for $n\ge 2$
\begin{equation}
\label{ref-A.5-48}
\begin{aligned}
\psi^n&=HPQ_2\psi^{\le n-1}I_n\\
Q^n_1&=pPQ_2\psi^{\le n-1}I_n
\end{aligned}
\end{equation}
Then
\begin{enumerate}
\item $Q_2\psi=\psi Q_1$. 
\item $(Q_1)^2=0$.
\item The functor $H^\ast(-)$ applied to  $\psi^1:H^\ast(V)\rightarrow V$ 
yields the identity on $H^\ast(V)$.
\item $Q_1^2:F^{c,2}(H^\ast(V))\rightarrow H^\ast(V)[1]$ coincides with $H^\ast(Q_2^2)$.
\end{enumerate}
\end{propositions}
\begin{proof}
(3)(4) are trivial so we concentrate on (1) and (2). 

(1) is equivalent to 
\[
PQ_2\psi I_n=P \psi Q_1 I_n
\]
for all $n\ge 1$. We prove this by induction on $n$, the case $n=1$
being clear.

We compute (with obvious notations)
\begin{align*}
PQ_2\psi I_n- P \psi Q_1 I_n&=Q^1_2\psi^nI_n+PQ_2\psi^{\le n-1}I_n
-\psi^1 Q^n_1 I_n-P\psi^{[2,n]}Q^{\le n-1}_1I_n\\
&=Q^1_2 HPQ_2\psi^{\le n-1}I_n+PQ_2\psi^{\le n-1}I_n
-\psi^1 pPQ_2\psi^{\le n-1}I_n-P\psi^{[2,n]}Q^{\le n-1}_1I_n\\
&=(\id+Q^1_2 H-\psi^1 p)PQ_2\psi^{\le n-1}I_n
-P\psi^{[2,n]}Q^{\le n-1}_1I_n\\
&=-HQ^1_2PQ_2\psi^{\le n-1}I_n
-P\psi^{[2,n]}Q^{\le n-1}_1I_n\\
&=-HPQ^1_2Q_2\psi^{\le n-1}I_n
-P\psi^{[2,n]}Q^{\le n-1}_1I_n\\
&=HPQ^{[2,n-1]}_2 Q_2^{\le n-1}\psi^{\le n-1}I_n-P\psi^{[2,n]}Q^{\le n-1}_1I_n\\
&=HPQ^{[2,n-1]}_2\psi^{\le n-1} Q_1^{\le n-1}I_n-P\psi^{[2,n]}Q^{\le n-1}_1I_n\qquad \text{(induction)}\\
&=P(HPQ^{[2,n-1]}_2\psi^{\le n-1}-\psi^{[2,n]}) Q_1^{\le n-1}I_n\\
&=0
\end{align*}
The argument for (2) is as in \cite{GH}. We include it for
completeness.  It is sufficient to prove
\[
P (Q_1)^2 I_{n+1}=0
\]
for all $n\ge 1$. Again $n=1$ is clear and we use induction for $n\ge 2$. 
Since $Q_1^1=0$ and $\psi^1$ is injective it
suffices to prove $\psi^1(Q_1^{\le n})^2I_{n+1}=\psi^{\le n}(Q_1^{\le
  n})^2I_{n+1}=0$. Since $Q_1^{\le n}$ maps $F^{c,n+1}(H^\ast(V))$ to
$F^{c,\le n}(H^\ast(V))$ we have
\begin{align*}
\psi^{\le
  n}Q_1^{\le n}Q_1^{\le n}I_{n+1}&=Q_2\psi^{\le
  n}Q_1^{\le n}I_{n+1}\qquad (\text{using (1)})\\
&=Q_2(\psi^{\le
  n}Q_1^{\le n}-Q_2\psi^{\le n})I_{n+1}\\
&=Q_2^1(\psi^{\le
  n}Q_1^{\le n}-Q_2\psi^{\le n})I_{n+1}
\end{align*}
where the last equality follows from the fact that $(\psi^{\le
  n}Q_1^{\le n}-Q_2\psi^{\le n})I_{n+1}$ takes values in $V$ (again using (1)).

We conclude
\[
\im \psi^1(Q_1^{\le n})^2 I_{n+1}\subset\im \psi^1\cap \im Q_2^1=0
\]
which finishes the proof.
\end{proof}
Combining \eqref{ref-A.4-46} with Proposition \ref{ref-A.2.1-47}
the recursion relations can be written as
\begin{multline}
\label{ref-A.6-49}
\psi^n(\phi\otimes v_1\otimes \cdots\otimes v_n)=\\
\sum_{\begin{smallmatrix}m\ge 2,\sigma\in S_n\\ \sum_{i=1}^m n_i=n\end{smallmatrix}}\pm
\frac{1}{m!\,n_1!\cdots n_m!} HQ^m_2((\sigma\phi)_{(0)}\otimes \psi^{n_1}\bigl((\sigma\phi)_{(1)}\otimes
  v_{\sigma^{-1}(1)}\otimes\cdots \otimes
  v_{\sigma^{-1}(n_1)}\bigr)\otimes\\\cdots \otimes
  \psi^{n_m}\bigl((\sigma\phi)_{(m)}\otimes
  v_{\sigma^{-1}(n_1+\cdots+n_{m-1}+1)}\otimes\cdots \otimes
  v_{\sigma^{-1}(n)}\bigr))
\end{multline}
\begin{multline}
\label{ref-A.7-50}
Q_1^n(\phi\otimes v_1\otimes \cdots\otimes v_n)=\\
\sum_{\begin{smallmatrix}m\ge 2,\sigma\in S_n\\ \sum_{i=1}^m n_i=n\end{smallmatrix}}\pm
\frac{1}{m!\,n_1!\cdots n_m!} pQ^m_2((\sigma\phi)_{(0)}\otimes \psi^{n_1}\bigl((\sigma\phi)_{(1)}\otimes
  v_{\sigma^{-1}(1)}\otimes\cdots \otimes
  v_{\sigma^{-1}(n_1)}\bigr)\otimes\\\cdots \otimes
  \psi^{n_m}\bigl((\sigma\phi)_{(m)}\otimes
  v_{\sigma^{-1}(n_1+\cdots+n_{m-1}+1)}\otimes\cdots \otimes
  v_{\sigma^{-1}(n)}\bigr))
\end{multline}

Let $\Tscr_n$ be the set
of planar rooted trees whose leaves are indexed from $1$ to $n$ and whose
internal vertices have at least two branches. 
Iterating \eqref{ref-A.6-49}\eqref{ref-A.7-50} it is clear that we will get
$\psi^n(\phi\otimes v_1\otimes\cdots\otimes v_n)$, $Q_1^n(\phi\otimes v_1\otimes\cdots\otimes v_n)$ as sums indexed by
elements of $S_n\times \Tscr_n$. 

To be more precise
 write for $v_1,\ldots, v_n\in V$
\[
\phi(v_1,\ldots,v_n)=Q_2^n(\phi\otimes v_1\otimes\cdots \otimes v_n)
\]
Then we get 
\begin{equation}
\label{ref-A.8-51}
\psi^n(\phi\otimes v_1\otimes \cdots\otimes v_n)=\\
\sum_{\begin{smallmatrix}\sigma\in S_n
\\ T\in \Tscr_n\end{smallmatrix}}\pm
\frac{1}{w_T}\psi_T(\sigma\phi\otimes v_{\sigma^{-1}(1)}\otimes\cdots\otimes  
v_{\sigma^{-1}(n)})
\end{equation}
\begin{equation}
\label{ref-A.9-52}
Q^n_1(\phi\otimes v_1\otimes \cdots\otimes v_n)=\\
\sum_{\begin{smallmatrix}\sigma\in S_n
\\ T\in \Tscr_n\end{smallmatrix}}\pm
\frac{1}{w_T}Q_{1,T}(\sigma\phi\otimes v_{\sigma^{-1}(1)}\otimes\cdots\otimes  
v_{\sigma^{-1}(n)})
\end{equation}
where $w_T$ is the product of the factorials of the internal vertices
and where $\psi_T$, $Q_{1,T}$ are inductively defined as follows: 
$\psi_\bullet=i$ (here $\bullet$ is the single vertex tree) and 
$$
\psi_{T}(\phi\otimes v_1\otimes\cdots\otimes v_n):=
H\phi_{(0)}\Big(
\psi_{T_1}\big(\phi_{(1)}\otimes v_1\otimes\cdots\otimes v_{|T_1|}\big),\dots,
\psi_{T_m}\big(\phi_{(m)}\otimes v_{1+\sum_{i=1}^{m-1}|T_i|}\otimes\cdots\otimes v_n\big)\Big)
$$
$$
Q_{1,T}(\phi\otimes v_1\otimes\cdots\otimes v_n):=
p\phi_{(0)}\Big(
\psi_{T_1}\big(\phi_{(1)}\otimes v_1\otimes\cdots\otimes v_{|T_1|}\big),\dots,
\psi_{T_m}\big(\phi_{(m)}\otimes v_{1+\sum_{i=1}^{m-1}|T_i|}\otimes\cdots\otimes v_n\big)\Big)
$$
if $T=B_+(T_1\cdots T_m)$ is obtained from $m$ trees by grafting them on a common 
new root. Here we have denoted by $|T|$ the number of leaves of a given tree $T$. 

We illustrate this with an example.
\begin{examples}
Let $T$ be the tree with corresponding parenthesized expression given
by $((1(23))4)$. Then we have 
\[
\psi_T(\phi\otimes v_1\otimes v_2\otimes v_3\otimes v_4)
=H\phi_{(0)}(H\phi_{(1)(0)}(i(v_1),H\phi_{(1)(2)(0)}(i(v_2),i(v_3))),i(v_4))
\]
where we have used $\Oscr(1)=k$ and $\psi^1=i$ and where
$\phi_{(0)}\otimes \phi_{(1)(0)}\otimes 1\otimes
\phi_{(1)(2)(0)}\otimes 1\otimes 1 \otimes 1$ is the image
of $\phi$ under the compositions
\begin{align*}
\Oscr(4)^\ast&\rightarrow \Oscr(2)^\ast\otimes \Oscr(3)^\ast\otimes \Oscr(1)^\ast\\
&\r\Oscr(2)^\ast\otimes (\Oscr(2)^\ast \otimes \Oscr(1)^\ast\otimes \Oscr(2)^\ast) \otimes \Oscr(1)^\ast\\
&\rightarrow \Oscr(2)^\ast\otimes \Oscr(2)^\ast \otimes \Oscr(1)^\ast\otimes (\Oscr(2)^\ast \otimes \Oscr(1)^\ast\otimes \Oscr(1)^\ast)\otimes \Oscr(1)^\ast
\end{align*}
Pictorially $\psi_T$ is given by
\[
\psfrag{v1}[][]{$v_1$}
\psfrag{v2}[][]{$v_2$}
\psfrag{v3}[][]{$v_3$}
\psfrag{v4}[][]{$v_4$}
\psfrag{i}[][]{$i$}
\psfrag{H}[][]{$H$}
\psfrag{a}[][]{$\phi_{(0)}$}
\psfrag{b}[][]{$\phi_{(1)(0)}$}
\psfrag{c}[][]{$\phi_{(1)(2)(0)}$}
\includegraphics{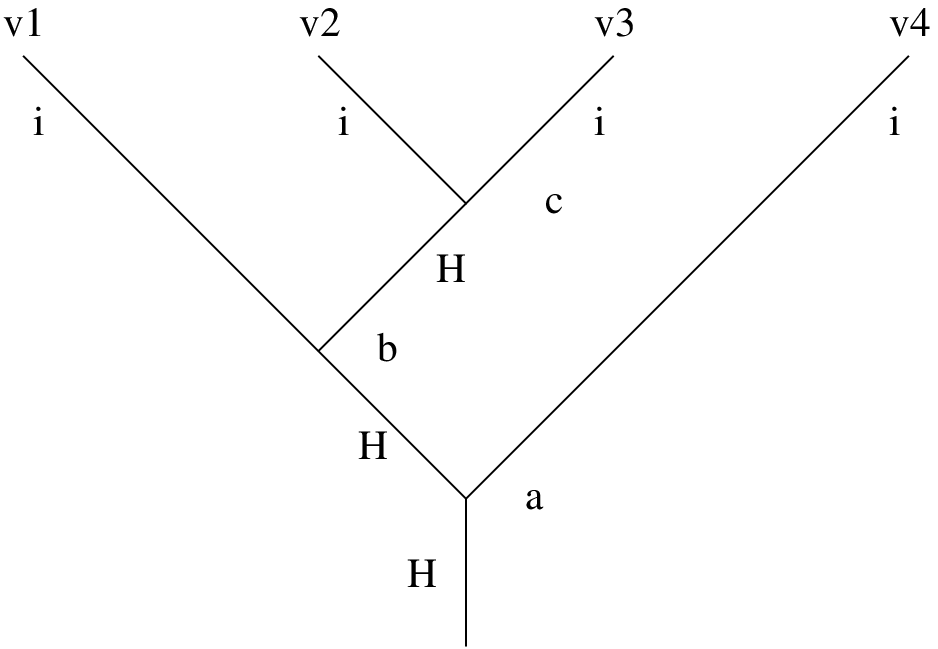}
\]
In this case $w_T=2!\times 2!\times 2!$. 
The pictorial representation for $Q_{1,T}$ is similar. The only difference
is that the root edge is decorated by $p$ instead of $H$. 
\end{examples}
Although it is not relevant for the discussion below we note that
\eqref{ref-A.8-51}\eqref{ref-A.9-52} can be more elegantly written in
terms of arbitrary trees instead of planar trees. To be more precise
let $\widetilde{\Tscr}_n$ be the set isomorphism classes of trees with
leaves given by $\{1,\ldots,n\}$. Then there is a map
\[
\gamma:S_n\times \Tscr_n\rightarrow \widetilde{\Tscr}_n
\]
which associates to $(\sigma, T)\in S_n\times \Tscr_n$ the tree obtained
from $T$ by replacing the leaf vertices by $1,\ldots,n$ by $\sigma^{-1}(1),
\ldots,\sigma^{-1}(n)$. 

It is clear that $\gamma$ is surjective but not injective. The
cardinalities of the fibers are precisely given by the numbers $w_T$
introduced above. 

If $\tilde{T}=\gamma(\sigma,T)$ then we put
\begin{align*}
\psi_{\tilde{T}}(\phi\otimes v_1\otimes\cdots\otimes  
v_n)&=\pm\psi_{T}(\sigma(\phi)\otimes v_{\sigma^{-1}(1)}\otimes
\cdots
\otimes
v_{\sigma^{-1}(n)})\\
Q_{1,\tilde{T}}(\phi\otimes v_1\otimes\cdots\otimes  
v_n)&=\pm Q_{1,T}(\sigma(\phi)\otimes v_{\sigma^{-1}(1)}\otimes
\cdots
\otimes
v_{\sigma^{-1}(n)})
\end{align*}
One may check that this is well defined using the operad axioms.
\eqref{ref-A.8-51}\eqref{ref-A.9-52} then become
\[
\psi^n(\phi\otimes v_1\otimes \cdots\otimes v_n)=\\
\sum_{\begin{smallmatrix}
\\ T\in \widetilde{\Tscr}_n\end{smallmatrix}}
\psi_T(\phi\otimes v_{1}\otimes\cdots\otimes  
v_{1})
\]
\[
Q^n_1(\phi\otimes v_1\otimes \cdots\otimes v_n)=\\
\sum_{\begin{smallmatrix}
\\ T\in \widetilde{\Tscr}_n\end{smallmatrix}}
Q_{1,T}(\phi\otimes v_{1}\otimes\cdots\otimes  
v_{n})
\]

\subsection{The construction of $\Psi'$}
We revert to the notations from the introduction. 
We prove the following result
\begin{propositions} \label{ref-A.3.1-53} There exists a $G_\infty$-structure $Q_1$ on
  $\frak{h}$ and  $G_\infty$-quasi-isomorphism
\[
\Psi':(\frak{h},Q_1)\rightarrow (\frak{g},Q_2)
\]
such that $\Psi'$ satisfies the obvious analogues of (P1-3);
$H^\ast(Q_2^{1,1})=\HKR\circ Q^{1,1}_1\circ \HKR^{-1}$ and
$H^\ast(Q^2_2)=\HKR\circ Q^2_1\circ \HKR^{-1}$; the components of $Q_1$ are
$\Aff(W)$-invariant poly-differential operators and in addition:
\begin{itemize}
\item[(R1)]
  $Q_1^{1,p_2,\ldots,p_n}(\gamma\underline{\alpha}_2\cdots\underline{\alpha}_n)=0$
  for $\gamma\in \frak{h}_0$
  and $\underline{\alpha}_i\in L^{c,p_i}(\frak{h})$, except when $n=2$ and
$p_2=1$.
\item[(R2)] $\Psi^{\prime 1,p_2,\ldots,p_n}(\gamma\underline{\alpha}_2\cdots\underline{\alpha}_n)=0$ for
  $\gamma\in \frak{h}_0$
  and $\underline{\alpha}_i\in L^{c,p_i}(\frak{h})$  and $n\geq2$.
\end{itemize}
Note that in contrast to the statement of (P5) in this case $\gamma$ is an arbitrary, not necessarily linear,
vector field.
\end{propositions}
\begin{proof}
Define $\Fscr^{m,n}(\frak{g})\subset \frak{g}$ as the vector space of
poly-differential operators of degree $\le n$ with
$m$ arguments. Clearly $\Fscr^{m,n}(\frak{g})$ is a finite free
$R$-module. A linear map $\frak{g}\times\cdots\times \frak{g}\r
\frak{g}$ is called a local $R$-poly-differential operator if its restriction
to any $\Fscr^{m_1,n_1}(\frak{g})\times \Fscr^{m_2,n_2}(\frak{g})\times \cdots$
is an $R$-poly-differential operator. Local poly-differential operators
are compatible with all operations we use below.

\medskip

It is easy to verify that structure maps in the $B_\infty$-structure
on $\frak{g}$ are $\Aff(W)$-invariant local $R$-poly-differential
operators.
Applying the morphism of DG-operads $G_\infty\rightarrow B_\infty$ we obtain 
 the same result for the $G_\infty$-structure on $\frak{g}$.

\medskip

If $V$ is a graded vector spaces then the cofree coalgebra cogenerated
by $V[1]$ for the Gerstenhaber algebra operad is equal to
$S^c(L^c(V)[1])$.  Hence if $V$ is a $G_\infty$-algebra represented by a
codifferential of degree one on $S(L(V)[1])$ then we may use the
\eqref{ref-A.8-51}\eqref{ref-A.9-52} to construct a minimal model for
$V$.

We apply this with $V=\frak{g}$ and we put
$i=\text{HKR}$. We claim that we may choose
$H$ and $p$ to be $\Aff(W)$-equivariant and $R$-linear such that in addition
we have $Hi=0$.

It is easy to write down a formula for $p$ but $H$ is another matter.
There exists explicit, but quite non-trivial, formulas for $H$ which
have the required properties \cite{Halbout1,dWL}. As the
explicit form of $H$ is not required for us we may also use the following
non-explicit argument.

As $\Aff(W)$-representations
\[
T_{\poly}(R)=R\otimes_k \wedge\overline{W}^\ast
\]
and 
\[
D_{\poly}(R)=R\otimes_k T_k(S\overline{W}^\ast)
\]
and $i=\operatorname{HKR}$ is clearly obtained by base extension from its
restriction
\begin{gather*}
i':\wedge \overline{W}^\ast\rightarrow T_k(S\overline{W}^\ast)
\end{gather*}
In particular $i'$ is still a quasi-isomorphism.
We choose a $\GL(W)$ invariant quasi-inverse 
\[
p': T_k(S\overline{W}^\ast)\r\wedge \overline{W}^\ast
\]
such that $p'i'=\Id$. This is possible since $\GL(W)$ is reductive. Then we choose a
$\GL(W)$-equivariant homotopy
\[
H':T_k(S\overline{W}^\ast)\rightarrow T_k(S\overline{W}^\ast)[-1]
\]
between the identity and $i'p'$ satisfying $H'i'=0$. Finally we let $p$, $H$ be the base extensions
of $p'$, $H'$.

Now we construct $Q_1$ and $\Psi'$ using \eqref{ref-A.8-51}
\eqref{ref-A.9-52}.  We obtain
(R1) and (R2) (which are {\it a priori} not part of the conclusions of Proposition
\ref{ref-A.3.1-53}) using the following two facts:
\begin{itemize}
\item $Q_2^{1,p_2,\ldots,p_n}(i(\gamma)\underline{\alpha}_2\cdots\underline{\alpha}_n)=0$
  for $\gamma\in \frak{h}_0$ and $\underline{\alpha}_i\in
  L^{c,p_i}(\frak{g})$, except when $n=2$ and $p_2=1$. 
This follows from Proposition \ref{ref-5.5-11} and Lemma \ref{ref-5.2-10}. 
\item $HQ_2^{1,1}(i(\gamma) i(\gamma'))=0$ for $\gamma\in \frak{h}_0$
  and $\gamma'\in \frak{h}$. This follows from the fact that
  $HQ_2^{1,1}(i(\gamma)
  i(\gamma'))=H[i(\gamma),i(\gamma')]=Hi([\gamma,\gamma'])=0$ where
the second equality follows from the fact that the HKR map is 
compatible with $[\delta,-]$ when $\delta$ is a vector field. \qed
\end{itemize}
\def\qed{}\end{proof}

\subsection{Sketch of Tamarkin's argument}
\label{ref-A.4-54}
We will first remind the reader of Tamarkin's argument (without
supplying all details) which yields a $G_\infty$-isomorphism
$\Psi'':(\frak{h},Q_0)\rightarrow (\frak{h},Q_1)$ such that $\Psi^{\prime
  \prime}=\Id_{\frak{h}}$, but satisfying a priori no additional
conditions. Then we will modify Tamarkin's argument to make $\Psi''$
satisfy stronger conditions.

\medskip

Put $F=S^c(L^c(\frak{h})[1])$. 
We first consider the ``deformation complex'' of $\frak{h}$
\begin{equation}
\label{ref-A.10-55}
(\coder(F,F)^\delta,[Q_0,-])= (\Hom(F,\frak{h}[1]),[Q_0,-])
\end{equation}
\begin{propositions}[Tamarkin \cite{Tamarkin}] \label{ref-A.4.1-56}
  The deformation complex of $\frak{h}$ is exact, except for a copy of $k$
  located in cohomological degree $-2$ which has as generating vector the
  composition $F\twoheadrightarrow k\hookrightarrow \frak{h}[1][-2]=S(W\oplus W^\ast[-1])$.
\end{propositions}
We sketch the proof below but for further reference we first 
state and prove the motivating corollary.
\begin{corollarys}
\label{ref-A.4.2-57}
There exists a $G_\infty$-isomorphism
$\Psi'':(\frak{h},Q_0)\rightarrow (\frak{h},Q_1)$ such that $\Psi^{\prime
  \prime 1}=\Id_{\frak{h}}$.
\end{corollarys}
\begin{proof}
We consider $F$ to be graded by 
$\frak{h}$-homogeneity. Below superscripts refer to this grading. 
Clearly $Q_0$ raises the grading by one\footnote{This depends crucially
on the fact that $Q_0$ is built only of binary operations. I.e. $\frak{h}$ it is not itself a strong homotopy algebra.} on $F$.  Hence $[Q_0,-]$ maps
$\Hom(F^n,\frak{h}[1])$ to $\Hom(F^{n+1},\frak{h}[1])$. 

In other words, we deduce that there is a complex of complexes
 with zero differential
\begin{equation}
\label{ref-A.11-58}
0\rightarrow k[2]\rightarrow \Hom(F^0,\frak{h}[1])\xrightarrow{[Q_0,-]} \Hom(F^1,\frak{h}[1])
\xrightarrow{[Q_0,-]}\Hom(F^2,\frak{h}[1])\xrightarrow{[Q_0,-]} 
\Hom(F^3,\frak{h}[1])\rightarrow \cdots
\end{equation}
Since the total complex is exact by Tamarkin's result, this is in fact an
exact sequence of complexes with zero differential. 

\medskip

Now we perform a standard computation.  Choose $n\ge 3$ minimal such
that $(Q_1-Q_0)|F^{n}\neq 0$. Then we have on $F^{\le n}$:
$0=[Q_1,Q_1]=[Q_0,Q_0]+2[Q_0,Q_1-Q_0]=2[Q_0,Q_1-Q_0]$. Hence by
\eqref{ref-A.11-58} we find $Z\in \Hom(F^{n-1},\frak{h}[1])$ such that
$[Q_0,H]=Q_1-Q_0$. Put $\theta=e^Z$. Then we find on $F^n$: $\theta
Q_1\theta^{-1}=Q_0-[Q_0,Z]+(Q_1-Q_0)=Q_0$.

Replacing $Q_1$ by $\theta Q_1\theta^{-1}$ and repeating this procedure
we eventually find the requested $G_\infty$-isomorphism. 
\end{proof}
\begin{proof}[Proof of Proposition \ref{ref-A.4.1-56} (sketch)]
The crucial
point is that $Q_0=Q^{\cup}_0+Q^{\Lie}$ where $Q^{\cup}_0$ and
$Q_0^{\Lie}$ are respectively obtained from the product and Lie
bracket on $\frak{h}$ and where $[Q_0^{\Lie},Q_0^{\cup}]=0$. 
\begin{enumerate}
\item The commutative algebra structure on $\frak{h}[-1]$ makes
  $\frak{h}[-1]$ into a $C_\infty$-algebra. This yields a
  codifferential $Q_0^{\text{cup}}$ on $L^c(\frak{h})$. Using Leibniz
rule the codifferential may be extended to a codifferential on $F$, also
denoted by $Q_0^{\text{cup}}$. 
\item The Lie algebra structure on $\frak{h}$ extends to a Lie algebra
structure on $L^c(\frak{h})$ such that $L^c(\frak{h})\rightarrow \frak{h}$ is
a Lie algebra morphism. For the explicit formula of the Lie bracket
(which is not entirely obvious) we refer to \cite[\S3.1]{Fresse}.

In particular $L^c(\frak{h})$ is an $L_\infty$-algebra and hence
there is a corresponding codifferential $Q_0^{\Lie}$ on $F$. 
\end{enumerate}
Using these explicit description we may compute the actions of 
$[Q_0^{\text{cup}},-]$ and $[Q_0^{\text{Lie}},-]$ on $\Hom(F,\frak{h}[1])$.
This is somewhat tricky to get right and we only list the results. We refer to Tamarkin's
paper for details. 
\begin{enumerate}
\item $[Q_0^{\Lie},-]$ is simply the Cartan-Eilenberg differential for the
$L^c(\frak{h})$ representation $\frak{h}$. 
\item To describe $[Q_0^{\cup},-]$ we put for clarity
  $A=\frak{h}[-1]=S(W\oplus W^\ast[-1])$ and consider $A$ as an
  associative algebra.
Then
\begin{align*}
\Hom(F,\frak{h}[1])&=\Hom_A(A\otimes_k F,\frak{h}[1])\\
&=\Hom_A(S^c_A(A\otimes_k L^c(\frak{h})[1]),\frak{h}[1])
\end{align*}
On $A\otimes_k L^c(\frak{h})=A\otimes_k L^c(A[1])$ we have the A-linear Harrison differential $d_{\Harr}$
which is obtained from the Hochschild differential on $A\otimes T^c(A[1])$
(using the surjection $A\otimes_k T^c(A[1])\rightarrow A\otimes_k L^c(A[1])$). 

We extend $d_{\Harr}$ to a differential on $S^c_A(A\otimes_k L^c(\frak{h})[1])$
which we also denote by $d_{\Harr}$. Then $[Q^{\cup}_0,-]$ is the
dual of $d_{\Harr}$. 
\end{enumerate}
We now relate the complex $(\Hom_A(S^c_A(A\otimes_k
L^c(\frak{h})[1]),\frak{h}[1]),[Q_0,-])$ to the complex computing the
Lichnerowisz-Poisson cohomology of $A$. This is explained
\cite[\S1.4.9]{Fresse} but we present a slightly different point of
view in terms of Lie algebroid cohomology.

If $\frak{l}$ is a Lie algebra acting on a commutative ring $S$ then
$S\otimes_k \frak{l}$ carries a natural structure of a Lie algebroid.
In our case the Lie algebra $\frak{h}$ acts on $A$ (via the Lie
bracket on $\frak{h}=A[-1]$). Hence $L^c(\frak{h})$ acts on $A$ via
the projection map $L^c(\frak{h})\rightarrow \frak{h}$ which is a Lie algebra
homomorphism. It follows $A\otimes L^c(\frak{h})$ is a Lie algebroid.
Using the explicit formula for the Lie bracket on $L^c(\frak{h})$ in
\cite[\S3.1]{Fresse} one checks that $(A\otimes_k
L^c(\frak{h}),d_{\text{Harr}})$ is in fact a DG-Lie algebroid. The
complex $(\Hom_A(S^c_A(A\otimes_k
L^c(\frak{h})[1]),\frak{h}[1]),[Q_0,-])$ is the complex computing
the cohomology of this Lie algebroid.

As $\frak{h}$ is a (shifted) Gerstenhaber algebra we find that
$\Omega_{A}[1]$ is a Lie algebroid with Lie bracket and anchor map
determined by $[df,dg]=d[f,g]$, $\rho(df)(g)=[f,g]$. The Lie algebroid
cohomology of $\Omega_A[1]$ is the (shifted) Lichnerowisz-Poisson
cohomology of $A$. 
Furthermore one has  morphisms of DG-Lie algebroids
\begin{equation}
\label{ref-A.12-59}
(A\otimes_k L^c(A[1]),d_{\text{Harr}})\rightarrow (A\otimes A[1],0)\xrightarrow{a\otimes b\mapsto adb} 
(\Omega_A[1],0)
\end{equation}
Dualizing we obtain a morphism of complexes
\begin{equation}
\label{ref-A.13-60}
\Hom_A(S^c_A(\Omega_A[2]),\frak{h}[1]),d_{\text{Poiss}})
\r
(\Hom_A(S^c_A(A\otimes_k L^c(\frak{h})[1]),\frak{h}[1]),[Q_0,-])
\end{equation}
where $d_{\Poiss}$ is the differential computing the (shifted) Lichnerowisz-Poisson
cohomology of $A$

We claim that \eqref{ref-A.13-60} is in fact a quasi-isomorphism. Using a spectral
sequence argument it is sufficient to prove that 
\[
\Hom_A(S^c_A(\Omega_A[2]),\frak{h}[1]),0)
\r
(\Hom_A(S^c_A(A\otimes_k L^c(\frak{h})[1]),\frak{h}[1]),[Q_0^{\cup},-])
\]
is a quasi-isomorphism. This follows from the fact that
\eqref{ref-A.12-59} is a homotopy equivalence. To prove this last statement
we note that it is well-known that \eqref{ref-A.12-59} is a quasi-isomorphism
(e.g.\ \cite[Lemma 2.5.10]{Kapra}). Then it suffices to observe that
$A\otimes L^c(A[1])$ is a colimit of finite extensions of shifts of $A$
and hence is homotopically projective. Thus a quasi-isomorphism 
with source $A\otimes L^c(A[1])$ is automatically a homotopy equivalence.

Using all this we are reduced to computing the (shifted) Lichnerowisz-Poisson
cohomology of $A$. Now the Gerstenhaber bracket is in fact symplectic
and it is well-known that Lichnerowisz-Poisson cohomology for a symplectic form is the
same as De Rham cohomology (see \cite{Bryl}). The fact that
the De Rham cohomology of $A$ is trivial then finishes the proof of 
Proposition \ref{ref-A.4.1-56}. We refer to \cite{Tamarkin} for more
details.
\end{proof}

\subsection{The differential deformation complex}

We define $\Diff(S^c(L^{c}(\frak{h})[1]),\frak{h}[1])$ as the graded subspace
of $\Hom(S^c(L^{c}(\frak{h})[1]),\frak{h}[1])$ such
that all the occurring maps $L^{c,m_1}(\frak{h})\otimes \cdots \otimes
L^{c,m_1}(\frak{h})\rightarrow \frak{h}$ are $R$-poly-differential operators.
\begin{lemmas}
$\Diff(S^c(L^{c}(\frak{h})[1]),\frak{h}[1])$ is closed under $[Q_0^{\cup},-]$
and $[Q_0^{\Lie},-]$.
\end{lemmas}
\begin{proof} Ultimately $[Q_0^{\cup},-]$ and $[Q_0^{\Lie},-]$ are derived
from the cupproduct and Lie bracket on $\frak{h}$. Since the latter
are $R$-poly-differential operators, one obtains that the same holds
for $[Q_0^{\cup},-]$ and $[Q_0^{\Lie},-]$.
\end{proof}
We have two aims in this section. The first one is to give a more
convenient expression for the differential deformation complex. See
\eqref{ref-A.15-63} and Lemma \ref{ref-A.5.2-62}.
The second aim is to show that the analogue of 
Proposition \ref{ref-A.4.1-56} holds for the differential deformation
complex. See \eqref{ref-A.16-64}.  These two results will be used
in the next section. 

\medskip

We have
\begin{align*}
\Diff(L^{c}(A[1]),A)&\hookrightarrow \Diff(T^{c}(A[1]),A)\\
&=T_A(\Diff(A,A)[-1])\\
&=D_{\poly}(A)[-1]
\end{align*}
Since $A=S(W\oplus W^\ast[-1])$ we have\footnote{It is easy to see
  that $A$-differential operators are the same as
  $R$-differential operators.} $\Diff(A,A)=A\otimes S(W^\ast\oplus
W[1])$ and thus
\begin{equation}
\label{ref-A.14-61}
T_A(\Diff(A,A)[-1])=A\otimes T(S(W^\ast\oplus W[1])[-1])
\end{equation}
\begin{lemmas} 
\label{ref-A.5.2-62} Denote the free Lie algebra generated by a graded vector
  space $V$ by $L(V)$. Then with the identification
  \eqref{ref-A.14-61} we have
\[
\Diff(L^{c}(A[1]),A)=A\otimes L(S(W^\ast\oplus W[1])[-1])
\]
\end{lemmas}
\begin{proof}
Since $L^c(A)$ is equal to $T^c(A)$ modulo shuffles, one 
quickly establishes that $\Diff(L^{c}(A[1]),A)$ is the set of primitive
elements for the coshuffle coproduct on $T_A(\Diff(A,A)[-1])$. 
This yields the desired result.
\end{proof}
Put $\Lscr(A)=\Diff(L^{c}(A[1]),A)$. Then we get
\begin{equation}
\label{ref-A.15-63}
\Diff(S^c(L^{c}(\frak{h})[1]),\frak{h}[1])=S_A(\Lscr(A)[-1])[2]
\end{equation}
The differential $[Q_0^{\cup},-]$ on $S_A(\Lscr(A)[-1])[2]$ is obtained
by extending the Harrison differential on $\Lscr(A)$ (obtained
from the Hochschild differential on $D_{\poly}(A)$). 

The HKR quasi-isomorphism
\[
S_A(\Hom_A(\Omega_A,A)[-1])\rightarrow D_{\poly}(A)[-1]
\]
restricts to a quasi-isomorphism
\[
\Hom_A(\Omega_A,A)[-1]\rightarrow \Lscr(A)
\]
so that we get a quasi-isomorphism
\[
(S_A(\Hom_A(\Omega_A,A)[-2]),0)\rightarrow (S_A(\Lscr(A)[-1]),[Q_0^{\cup},-])
\]
and one checks that this is compatible with $Q_0^{\Lie}$ using the fact
that it is the corestriction of (a shifted version of) \eqref{ref-A.13-60}.

So ultimately using \eqref{ref-A.13-60} we obtain a quasi-isomorphism
\[
(S_A(\Hom_A(\Omega_A[2],A))[2],d_{\Poiss})\rightarrow (S_A(\Lscr(A)[-1])[2],[Q_0,-])=
 (\Diff(S^c(L^{c}(\frak{h})[1]),\frak{h}[1]),[Q_0,-])
\]
and since the cohomology is of the left hand side is $k[2]$ (see \S\ref{ref-A.4-54}) one obtains that the analogue of Proposition \ref{ref-A.4.1-56} holds
for the differential deformation complex. 
\begin{equation}
\label{ref-A.16-64}
(\Diff(S^c(L^{c}(\frak{h})[1]),\frak{h}[1]),[Q_0,-])\cong k[2]
\end{equation}
\begin{remarks}   Being a symmetric algebra the ring $R$ has a
  natural ascending filtration. So $\Diff(F,\frak{h}[1])$ has a corresponding
  descending filtration (we view $\frak{h}[1]$ as not filtered) and the completion
for this filtration is $\Hom(F,\frak{h}[1])$. It seems not unlikely that (with some
more work) this
observation could be used to deduce
  \eqref{ref-A.16-64} from Proposition \ref{ref-A.4.1-56}.
\end{remarks}

\subsection{The construction of  $\Psi''$}
\label{ref-A.6-65}

Consider
\[
\bar{\frak{d}}=\{Q\in
\coder(S^c(L^c(\frak{h})[1]),S^c(L^c(\frak{h})[1]))^{\delta}
\mid
\forall \gamma\in \afff(W):Q^{1,p_2,\ldots,p_n}(\gamma,\ldots)=0\}
\]
Then $\bar{\frak{d}}$ is clearly a Lie subalgebra of
$\coder(S^c(L^c(\frak{h})[1]),S^c(L^c(\frak{h})[1]))^{\delta}$ and
hence so is $\frak{d}\overset{\text{def}}{=}\bar{\frak{d}}^{\Aff(W)}$.

We have $Q^{\text{cup}}_0\in \frak{d}$ and hence $\bar{\frak{d}}$ and $\frak{d}$ are closed
under $[Q^{\text{cup}}_0,-]$. On the other hand $Q^{\Lie}_0\not\in \bar{\frak{d}}$.
Nonetheless one checks that $\frak{d}$ is closed under $[Q^{\Lie}_0,-]$ (but
not $\bar{\frak{d}}$).

In this way we obtain the ``affine equivariant deformation complex''
\[
(\frak{d},[Q_0,-])=\Hom(S^c(L^c(\frak{h})/\afff(W)[1]),\frak{h}[1])^{\Aff(W)},[Q_0,-])
\]
This complex can be defined in yet another way. For simplicity write
$d=[Q_0,-]$.

For $\gamma\in \afff(W)$ define the endomorphism $i_\gamma$ of degree $-1$
of $\Hom(S^c(L^c(\frak{h})[1]),\frak{h}[1])$ which sends $Q$ to
$Q(\gamma,\ldots)$. The one checks that $\gamma\in \afff(W)$ acts by
$L_\gamma\overset{\text{def}}{=}di_\gamma+i_\gamma d$ and hence
\[
\frak{d}=\{Q\in \Hom(S^c(L^c(\frak{h})[1]),\frak{h}[1]): \forall 
\gamma\in \afff(W)\mid i_\gamma Q=L_\gamma Q=0\}
\]
In this way we see the connection with equivariant cohomology. 

In order to construct a $\Psi''$ satisfying (P3)(P5) we would have to
analyze the cohomology of $\frak{d}$. However we want to satisfy also
(P1). Therefore we consider the subcomplex
\begin{equation}
\label{ref-A.17-66}
D=(\Diff(S^c(L^c(\frak{h})/\afff(W)[1]),\frak{h}[1])^{\Aff(W)},[Q_0,-])
\end{equation}
defined as usual by requiring that all occurring maps are $R$-poly-differential
operators. The cohomology of $D$ is given by Theorem \ref{ref-A.1.2-42}. We
will prove that theorem after we have proved the next corollary.  
\begin{corollarys} (to Theorem
  \ref{ref-A.1.2-42}) \label{ref-A.6.1-67} Assume that $Q_1$ was
  chosen as in Proposition \ref{ref-A.3.1-53}. Then there exists a
  $G_\infty$-isomorphism $\Psi'':(\frak{h},Q_0)\rightarrow (\frak{h},Q_1)$ such
  that $\Psi^{\prime \prime 1}=\Id_{\frak{h}}$ and such that the
  obvious analogues of (P1)(P3)(P4)(P5) hold for $\Psi''$.
\end{corollarys}
\begin{proof} Of course we proceed as in Corollary \ref{ref-A.4.2-57}. 
The grading by $\frak{h}$-homogeneity we used on $\Hom(F,\frak{h}[1])$ induces
a grading on $D$. Then we have a complex of complexes with zero differential
\begin{equation}
\label{ref-A.18-68}
0\rightarrow D^0\xrightarrow{[Q_0,-]} D^1\xrightarrow{[Q_0,-]} D^2\xrightarrow{[Q_0,-]}
\end{equation}
Then by (R1) we have that $Q_1{\mid} F^{n}\in D^n$ for $n\ge 3$. Since
$Q_0{\mid} F^{n}=0$ for $n\ge 3$ and $Q_0{\mid} F^n=Q_1{\mid} F^n$ for
$n\le 2$ we find $Q_1-Q_0\in D$.

Choose $n\ge 3$ minimal such that $Q_1-Q_0{\mid} F^n\neq 0$.  As in
the proof of Corollary \ref{ref-A.4.2-57} we find $[Q_0,Q_1-Q_0]=0$
on $F^{n+1}$. Now the point is that $Q_1-Q_0$ has \emph{cohomological}
degree one (not Hochschild degree one). Since by Theorem
\ref{ref-A.1.2-42} the total complex of \eqref{ref-A.18-68} has no
cohomology in odd degree we find $Z\in D^{n-1}$ such that
$[Q_0,Z]=Q_1-Q_0|F^n$.

We put $Q'_1=
e^ZQ_1e^{-Z}=Q_1+\sum [\cdots [Z,Q_1]\cdots]$. Now $[Z,Q_1]=-[Q_0,Z]+[Z,Q_1-Q_0]\in D$. Hence $Q'_1\in D$. Replacing $Q'_1$ and iterating we find a $\Psi''$
satisfying (P1)(P3)(P5). 

We claim that (P4) is  satisfied. We first observe that
for $n>2$, (P4) is automatic for degree reasons. For $n=2$ we observe that
$\Psi^{\prime\prime 1,1}$ yields an affine invariant differential operator 
\begin{equation}
\label{ref-A.19-69}
\bigwedge^2 T_{\poly}(R)_0\rightarrow R
\end{equation}
Any such affine invariant differential operator $P$ is of the form 
$$
P=f^{i,j}_{I,J}(\underline{x})\partial_{\xi_i}\partial^I_{\underline{x}}\wedge\partial_{\xi_j}\partial^J_{\underline{x}}\,,
$$
where $\underline{x}=(x_1,\dots,x_n)$ is a coordinate system on $W^*$
and $\xi_i$'s are the corresponding odd coordinates on $W[1]$; and
$\partial^I_{\underline{x}}=\prod_{t=1}^k\partial_{x_{i_t}}$ for
$I=(i_1,\dots,i_k)$.

Translation invariance implies that $f^{i,j}_{I,J}$ is actually a constant polynomial, and $GL(W)$-invariance 
imposes that $P$ can only be (up to a scalar factor) $\partial_{\xi_i}\wedge(\partial_{\xi_j}\partial_{x_i}\partial_{x_j})$, which does not happen to satisfy (P5). 
The corollary is therefore proved.
\end{proof}
\begin{proof}[Proof of Theorem \ref{ref-A.1.2-42}]
Put 
\[
\Diff_0(A,A)=\{D\in \Diff(A,A)\mid D{\mid}\afff(W)=0\}
\]
\begin{sublemma}
There is a a split exact sequence of $(A,\Aff(W))$-modules
\[
0\r\Diff_0(A,A)\rightarrow \Diff(A,A)\rightarrow A\otimes \afff(W)^\ast\rightarrow 0
\]
where the rightmost non-trivial map is given by $D\mapsto D\mid \afff(W)$. 
\end{sublemma}
\begin{proof} This sequence is obviously left exact. To prove the claims
it is sufficient to construct the splitting of the right most non-trivial map.

Choose a basis $(t_i)_{i=1}^d$ for $W$. Then $A$ is the graded
commutative algebra generated by $t_i$ and
$\partial_i=\partial/\partial t_i$. Denote the partial derivatives on
$A$ with respect to $t_i$ and $\partial_i$ by $\partial_{t_i}$ and
$\partial_{\partial_i}$.

A basis of $\afff(W)\subset A$ is given by $(\partial_i)_i$ and
$(t_i\partial_j)_{ij}$. Let $E=\sum_i t_i\partial_{t_i}\in \Diff(A,A)$. 
We send the element of $A\otimes \afff(W)^\ast$, given by
$\partial_i\mapsto a_i$, $t_i\partial_j\mapsto a_{ij}$ to 
$\sum_i a_i(1-E)\partial_{\partial_i}+\sum_{ij}a_{ij}\partial_{t_i}\partial_{\partial_j}\in \Diff(A,A)$. 

This is an $A$-linear splitting and one checks that it is indeed
$\Aff(W)$ equivariant. 
\end{proof}
If we  put
\[
\Lscr_0(A)=\Diff(L^{c}(A[1])/\afff(W),A)
\]
we obtain an exact sequence
\begin{equation}
\label{ref-A.20-70}
0\rightarrow \Lscr_0(A)\rightarrow \Lscr(A)\rightarrow A\otimes \afff(W)^\ast\rightarrow 0
\end{equation}
still split as $(A,\Aff(W))$-modules. Similarly as in \eqref{ref-A.15-63} we get
\[
\Diff(S^c(L^c(\frak{h})/\afff(W)[1]),\frak{h}[1])^{\Aff(W)}=
(S_A(\Lscr_0(A)[-1]))^{\Aff(W)}[2]
\]
\begin{sublemma} The inclusion $S_A(\Lscr_0(A)[-1]))\hookrightarrow 
S_A(\Lscr(A)[-1])$ extends to a quasi-isomorphism
\begin{equation}
\label{ref-A.21-71}
(S_A(\Lscr_0(A)[-1]))^{\Aff(W)}\r
((S_A(\Lscr(A)[-1]))\otimes_k S(\afff(W)^\ast[-2]))^{\Aff(W)}
\end{equation}
where the differential on the righthand side is given by
\begin{equation}
\label{ref-A.22-72}
d\otimes 1+\sum_i i_{e_j}\otimes e_j^\ast
\end{equation}
with $(e_j)_j$ an arbitrary basis for $\afff(W)$. 
\end{sublemma}
\begin{proof} This is a standard observation in equivariant
  cohomology. One first check that the square of \eqref{ref-A.22-72} is indeed
zero and that \eqref{ref-A.21-71} is indeed a morphism of complexes.

To prove that it is a quasi-isomorphism we use an appropriate spectral
sequence argument to reduce to proving that 
\[
((S_A(\Lscr_0(A)[-1]))^{\Aff(W)},0)\rightarrow ((S_A(\Lscr(A)[-1])\otimes_k S(\afff(W)^\ast[-2]))^{\Aff(W)},\sum_j i_{e_j}\otimes e_j^\ast)
\]
is a quasi-isomorphism. Before taking $\Aff(W)$-invariants this is
obviously a quasi-isomorphism since it is simply a kind of Koszul
resolution. But since \eqref{ref-A.20-70} is split is actually an
$\Aff(W)$-equivariant homotopy equivalence. Hence it remains a
homotopy equivalence after taking invariants.
\end{proof}
Let $\Tscr(W)\subset \Aff(W)$ be the translation group and let $\frak{t}(W)$ be
its corresponding Lie algebra. 
\begin{sublemma} Assume that $V$ is a rational $\frak{t}(W)$-representation. Then
$V\otimes_k R$ is injective in the category of all $\frak{t}(W)$-representations
(not just rational ones). In particular
\[
\Ext^i_{\frak{t}(W)}(k,V\otimes R)=
\begin{cases}
(V\otimes R)^{\Tscr(W)}&\text{if $i=0$}\\
0&\text{otherwise}
\end{cases}
\]
\end{sublemma}
\begin{proof} Let $U(\frak{t}(W))$ be the enveloping algebra of
  $\frak{t}(W)$. The category of $\frak{t}(W)$-representations is nothing
  but the category of $U(\frak{t}(W))$-modules. We have $R=SW$ and $U(\frak{t}(W))=SW^\ast$ and the action of $U(\frak{t}(W))$ on $SW$ is given by contraction. It
is then well-known that $SW$ is an injective $U(\frak{t}(W))$-module. 

Since $U(\frak{t}(W))$ is noetherian, a direct limit of injectives is injective,
so we may assume that $V$ is finite dimensional. Then we have
\[
\Hom_{U(\frak{t}(W))}(-,V\otimes R)=\Hom_{U(\frak{t}(W))}(V^\ast\otimes-,R)
\]
which is an exact functor. So we are done. 
\end{proof}
We have a two-step $\Aff(W)$-invariant filtration on $\afff(W)^\ast$ given by
\[
0\rightarrow \frak{gl}(W)^\ast\rightarrow \afff(W)^\ast \rightarrow \frak{t}(W)^\ast\rightarrow 0
\]
which induces a filtration by ideals on 
$S(\afff(W)^\ast[-2])$ and hence a filtration on
$
S_A(\Lscr(A)[-1])\otimes_k S(\afff(W)^\ast[-2])
$
with associated graded given by 
\[
S_A(\Lscr(A)[-1])\otimes_k S(\frak{t}(W)^\ast[-2])\otimes_k S(\frak{gl}(W)^\ast[-2])
\]
Using the second sublemma we find that this is compatible with taking 
$\Aff(W)$-invariants.  Thus we obtain
\begin{equation}
\label{ref-A.23-73}
\begin{aligned}
(\gr (S_A(\Lscr(A)[-1])&\otimes_k S(\afff(W)^\ast[-2])))^{\Aff(W)}\\
&=
(S_A(\Lscr(A)[-1])\otimes_k S(\frak{t}(W)^\ast[-2])\otimes_k S(\frak{gl}(W)^\ast[-2]))^{\Aff(W)}\\
&=((S_A(\Lscr(A)[-1])\otimes_k S(\frak{t}(W)^\ast[-2]))^{\Tscr(W)}\otimes_k S(\frak{gl}(W)^\ast[-2]))^{\GL(W)}
\end{aligned}
\end{equation}
We now claim that the only cohomology of $((S_A(\Lscr(A)[-1])\otimes_k
S(\frak{t}(W)^\ast[-2]))^{\Tscr(W)}$ is a copy of $k$ in degree zero.
To this end we use the following sublemma.
\begin{sublemma}
The inclusion
\[
((S_A(\Lscr(A)[-1])\otimes_k S(\frak{t}(W)^\ast[-2]))^{\Tscr(W)}\hookrightarrow
S_A(\Lscr(A)[-1])\otimes_k S(\frak{t}(W)^\ast[-2])
\]
extends to a quasi-isomorphism
\begin{multline}
\label{ref-A.24-74}
(S_A(\Lscr(A)[-1])\otimes_k S(\frak{t}(W)^\ast[-2]))^{\Tscr(W)}
\\ \rightarrow S_A(\Lscr(A)[-1])\otimes_k S(\frak{t}(W)^\ast[-2])\otimes S(\frak{t}(W)^\ast[-1])
\end{multline}
where the differential on the righthand side is given 
\begin{equation}
\label{ref-A.25-75}
d\otimes 1\otimes 1+\sum_{j} i_{e_j}\otimes e_j^\ast\otimes 1+
\sum_{k} L_{e_k}\otimes 1\otimes e_k^\ast
-\sum_{l}1\otimes e_l^\ast\otimes \partial_{e_l}
\end{equation}
for an arbitrary basis $(e_i)_i$ of $\frak{t}(W)$. 
\end{sublemma}
Note that the complex \eqref{ref-A.25-75} is the (unrestricted) BRST
model for equivariant $\Tscr(W)$-cohomology
\cite[\S3]{Kalkman}. Normally we would need to take a certain
subcomplex to get the correct result but in this case this is not
necessary because of the $\frak{t}(W)$ injectivity of $R$ which is
essentially a manifestation of the contractibility of $\Tscr(W)$.
\begin{proof} One first checks that \eqref{ref-A.25-75} has square zero. Then
one filters \eqref{ref-A.24-74} according to $\frak{t}(W)^\ast$ homogeneity
in $S(\frak{t}(W)^\ast[-1])$. Then we have to show that
\begin{multline}
\label{ref-A.26-76}
((S_A(\Lscr(A)[-1])\otimes_k S(\frak{t}(W)^\ast[-2]))^{\Tscr(W)},0)
\\\rightarrow (S_A(\Lscr(A)[-1])\otimes_k S(\frak{t}(W)^\ast[-2])\otimes S(\frak{t}(W)^\ast[-1]),\sum_{k} L_{e_k}\otimes 1\otimes e_k^\ast)
\end{multline}
is a quasi-isomorphism. 

Now the righthand side of \eqref{ref-A.26-76} computes 
\[
\Ext^\ast_{\frak{t}}(k,(S_A(\Lscr(A)[-1])\otimes_k S(\frak{t}(W)^\ast[-2])))
\]
To see this replace $\frak{t}(W)$ by $U(\frak{t}(W))=SW^\ast$ and then replace
$k$ by its Koszul resolution over $SW^\ast$. 

Since $S_A^n(\Lscr(A))\cong V_n\otimes R$
for  suitable rational $\frak{t}(W)$ representations $V_n$ we may conclude
by the second sublemma above. 
\end{proof}
To finish the computation of the cohomology of $S_A(\Lscr(A)[-1])\otimes_k
S(\frak{t}(W)^\ast[-2]))^{\Tscr(W)}$ we observe that the inclusion $k\subset S_A(\Lscr(A)[-1]$ extends to a morphism of complexes
\begin{multline*}
k\otimes_k S(\frak{t}(W)^\ast[-2])\otimes S(\frak{t}(W)^\ast[-1])
\r
S_A(\Lscr(A)[-1])\otimes_k S(\frak{t}(W)^\ast[-2])\otimes S(\frak{t}(W)^\ast[-1])
\end{multline*}
where the righthand side is as above and the lefthand side has differential
$-\sum_{l}1\otimes e_l^\ast\otimes \partial_{e_l}$. We claim that this is
again a quasi-isomorphism. Considering an appropriate filtration 
this follows from the fact that $k\rightarrow \Lscr(\Ascr)[-1]$ is a quasi-isomorphism
by \eqref{ref-A.16-64}. 

So it now remains to compute the cohomology of
$(S(\frak{t}(W)^\ast[-2])\otimes S(\frak{t}(W)^\ast[-1]), \sum_l
e_l^\ast\otimes \partial_{e_l})$, but this is just an ordinary Koszul
complex, so its cohomology is $k$.

\medskip

We now return to the spectral sequence derived from \eqref{ref-A.23-73}. Its first page becomes by the above discussion
\[
k\otimes_k S(\frak{gl}(W)^\ast[-2]))^{\GL(W)}
\]
and after that it degenerates. This finishes the proof (taking into account
that $\frak{gl}(W)^\ast\cong \frak{gl}(W)$ as $\GL(W)$-representations). 
\end{proof}
The proof of Theorem now follows easily by putting $\Psi=\Psi'\Psi''$ and
combining Proposition \ref{ref-A.3.1-53} with Corollary \ref{ref-A.6.1-67}.
\subsection{Kontsevich graphs}
\label{ref-A.7-77}
In this section we describe in more detail the $R$-poly-differential
operators between the $R$-modules $ T_{\poly}(R)$ and $D_{\poly}(R) $
which are equivariant under the affine group. This material is
well-known to experts. See for example \cite{GamHal,Shoikhet}.

\medskip

An ``admissible'' graph (or ``Kontsevich graph'') is an oriented graph with
the following properties. 
\begin{enumerate}
\item There are $t$ vertices of the ``first type''  labeled by
$1,\ldots,t$. Vertex $v$ has $n_v$ outgoing edges.
\item There are $n$ vertices of the ``second type'', labeled by $1,\ldots,n$,
with no outgoing edges.
\end{enumerate}
We write $\Gamma_i$ for the vertices of type $i$ in an admissible
graph $\Gamma$.

\medskip

Fix a basis $(t_i)_i$ for $W$ and write $\partial_i=\partial/\partial
t_i$.
For $a_{v}=\sum_{s_1,\cdots, s_{n_v}}a_v^{s_1\cdots
    s_{n_v}}\partial_{s_1}\cdots \partial_{s_{n_v}} \in
  T^{n_v}_{\poly}(R)$  and $f_1,\ldots,f_n\in R$ we define
\[
\Uscr_\Gamma(a_1,\ldots,a_t)(f_1,\ldots,f_n)=
\prod_{%
\begin{smallmatrix}
v\in \Gamma_1\\
\text{In}(v)={r_1,\ldots,r_d}\\
\text{Out}(v)={s_1,\ldots,s_{n_v}}
\end{smallmatrix}
}
 \partial_{{r_1}}\cdots \partial_{{r_d}} a^{s_1\cdots s_{n_v}}
\prod_{
\begin{smallmatrix}
v\in \Gamma_2\\
\text{In}(v)={u_1,\ldots,u_e}
\end{smallmatrix}
}
\partial_{{u_1}}\cdots \partial_{{u_e}}f_v
\]
Then clearly $\Uscr_\Gamma(a_1,\ldots,a_t)$ is an element of
$D_{\poly}^{n}(R)$ and hence $\Uscr_\Gamma$ defines a map
\[
T^{n_1}_{\poly}(R)\times\cdots\times T^{n_t}_{\poly}(R)\rightarrow D^{n}_{\poly}(R)
\]
which is an $\Aff(W)$ invariant poly-differential operator. 
\begin{propositions}\label{ref-A.7.1-78} An affine invariant poly-differential operators
  $T_{\poly}(R)\rightarrow D_{\poly}(R)$ is a linear combination of operators
  $\Uscr_\Gamma$ where $\Gamma$ runs through the admissible graphs.
\end{propositions}
\begin{proof}
We have
as $\Aff(W)$ representations
\[
T^n_{\poly}(R)=R\otimes_k \wedge^n\overline{W}^\ast
\]
and 
\[
D^n_{\poly}(R)=R\otimes_k (S\overline{W}^\ast)^{\otimes n}
\]
Thus for an arbitrary $R$-module the $R$-poly-differential operators from
$T^{n_1}_{\poly}(R)\times\cdots \times T^{n_t}_{\poly}(R)$ to $N$ are
given by
\[
\wedge^{n_1}\overline{W}\otimes_k\cdots\otimes_k
  \wedge^{n_t}\overline{W} \otimes_k (S\overline{W}^\ast)^{\otimes t}\otimes N
\]
Hence the $\Aff(W)$-invariant $R$-poly-differential operators
from
$T^{n_1}_{\poly}(R)\times\cdots \times T^{n_t}_{\poly}(R)$ to $D^n_{\poly}(R)$
are given by
\begin{align*}
  M&=(\wedge^{n_1}\overline{W}\otimes_k\cdots\otimes_k
  \wedge^{n_t}\overline{W} \otimes_k (S\overline{W}^\ast)^{\otimes t}
  \otimes_k(S\overline{W}^\ast)^{\otimes n} \otimes_k R)^{\Aff(W)}\\
&=(\wedge^{n_1}\overline{W}\otimes_k\cdots\otimes_k
  \wedge^{n_t}\overline{W} \otimes_k (S\overline{W}^\ast)^{\otimes t}
  \otimes_k(S\overline{W}^\ast)^{\otimes n})^{\GL(W)}
\end{align*}
If follows from Schur-Weyl duality that an $M$ is spanned by elements
$m_\Gamma$ associated to ``admissible'' graphs $\Gamma$ where
$m_\Gamma$ is defined as follows.  Write $dt_i$ for the element $t_i$
considered as an element of $\overline{W}$. Then
\[
m_\Gamma=\prod_{\begin{smallmatrix}v\in \Gamma_1\\
\text{Out}(v)={s_1,\ldots,s_{n_v}}\end{smallmatrix}}
dt_{s_1}\cdots dt_{s_{n_v}}
\otimes 
\prod_{\begin{smallmatrix}v\in \Gamma_1\\
\text{In}(v)={r_1,\ldots,r_d}
\end{smallmatrix}}
\partial_{{r_1}}\cdots \partial_{{r_d}}
\otimes 
\prod_{%
\begin{smallmatrix}
v\in \Gamma_2\\
\text{In}(v)={u_1,\ldots,u_e}
\end{smallmatrix}
}
\partial_{{u_1}}\cdots \partial_{{u_e}}
\]

A moment inspection reveal that $m_\Gamma$ corresponds to $\Uscr_\Gamma$.
\end{proof}
\begin{remarks} Since the Euler operator is invariant under
$\GL(W)$ but not under $\Aff(W)$, the conclusion of 
Proposition \eqref{ref-A.7.1-78} does not hold if we only demand that the
poly-differential operators are invariant under $\GL(W)$. 
\end{remarks}


\def\cprime{$'$} \def\cprime{$'$} \def\cprime{$'$}
\ifx\undefined\bysame
\newcommand{\bysame}{\leavevmode\hbox to3em{\hrulefill}\,}
\fi

\end{document}